\documentclass[showpacs,showkeys,amssymb,pra,notitlepage]{revtex4-1}

\usepackage{graphicx}

\usepackage{amsmath}

\usepackage{fancyvrb}
\usepackage[pdfauthor={Richard J. Mathar},colorlinks=true,bookmarksopen=true,pdfpagelayout=OneColumn]{hyperref}

\begin{document}

\title[Orthogonal Basis Over the Unit Circle]{Orthogonal Basis Function Over the Unit Circle with the Minimax Property}

\author{Richard J. Mathar}
\homepage{http://www.mpia.de/~mathar}
\pacs{42.15.Fr, 95.75.Pq}
\affiliation{Max-Planck Institute of Astronomy, K\"onigstuhl 17, 69117 Heidelberg, Germany}

\date{\today}
\keywords{Circle Polynomial, Orthogonal Basis, Chebyshev}

\begin{abstract}
We construct an orthogonal basis of functions defined over the unit circle as the product
of the common sinusoidal functions of the azimuth angle  by radial functions which
are essentially sines of a polynomials of the radial distance to the origin.
The main impetus of this approach is to generate basis
functions where the minima and maxima along both coordinates, the azimuth and the distance $r$ to the center,
have the same amplitude, akin to the Chebyshev polynomial basis of the one-dimensional unit interval.
The construction is based on numerical evaluation of the overlap integrals, which have the format of generalized Fresnel integrals.
\end{abstract}

\maketitle

\section{Aim and Scope} 

The celebrated Zernike polynomials define an orthogonal system of functions that
are products $M_m(\varphi)R_n^m(r)$ of azimuthal functions $M_m$ over $0\le \varphi\le 2\pi$
by radial polynomials $R_n^m(r)$ depending on the distance to the circle center, $0\le r\le 1$ \cite{LakshmJMO58}.
Some weakness of that system is that the polynomials
start to oscillate with increasing amplitude near $r\to 1$ for large $n$, such that the coefficients
of the expansion of some function of interest over that basis do not convey well
how important these functions are compared to the amplitudes contributed by
the other bases concentrated nearer to the origin.

We shall replace these unevenly wiggly radial polynomials by functions $Q_{n,m}(r)$
that have the minimax property of the Chebyshev Polynomials of the first kind \cite{FraserJACM12},
which means by functions where successive minima and maxima along $r$ have equal absolute
value. The difference to the Chebyshev Polynomials is of course the weight
$r$ in the orthogonality relation of the overlap integrals:
\begin{equation}
\int_0^1 r Q_{n,m}(r)Q_{n',m}(r)dr=\delta_{n,n'}.
\label{eq.ortho}
\end{equation}
As the azimuthal function are likewise orthogonal,
\begin{equation}
\int_0^{2\pi} M_{m}(\varphi)M_{m'}(\varphi)d\varphi=\delta_{m,n'},
\end{equation}
\begin{equation}
M_{m}(\varphi) = \sqrt{\frac{\epsilon_m}{2\pi}}\times\left\{
\begin{array}{ll}
\cos(m\varphi), & m\, \mathrm{even};\\
\sin(m\varphi), & m\, \mathrm{odd},\\
\end{array}
\right.
\end{equation}
(with $\epsilon_0= 1$ and $\epsilon_{|m|\ge 1}=2$) the orthogonality of the $Q$-basis is only concerned with orthogonality with
respect to the index $n$.

Among various strategies to approach the minimax property, we select one possible route;
we abandon the idea of expressing $Q$ as polynomials---giving up
some of the associated nice analytical results of the Zernike basis \cite{TangoApplPhys13,MatharSAJ179,JanssenJEOS6}
or other orthogonal bases \cite{ShengJOSAA11}---
and enforce the minimax property by using sines and cosines also along the radial coordinate.

\section{Design of Circle Functions}

For odd $m$ we construct radial functions
\begin{equation}
Q_{n,m}(r) = N_{n,m} \sin \theta_{n,m}(r),\quad n=m,m+2,m+4,\ldots,\quad m\,\mathrm{odd}
\end{equation}
and for even $m$ we construct radial functions
\begin{equation}
Q_{n,m}(r) = \left\{ 
\begin{array}{ll}
N_{n,m} \sin \theta_{n,m}(r),\quad n=m,m+2,m+4,\ldots, \quad m>0,\quad m\, \mathrm{even} \\
N_{n,0} \cos \theta_{n,0}(r),\quad n=0,2,4,\ldots , \quad m=0\\
\end{array}
\right.
\end{equation}
The cosine terms of $m=0$ are chosen to yield nonzero values at the origin, representing
forms of astigmatism.
The parities are kept in accordance with the Zernike polynomials:
\begin{equation}
\theta_{n,m}(-r) = (-)^m\theta_{n,m}(r),
\end{equation}
and therefore
\begin{equation}
Q_{n,m}(-r) = (-)^mQ_{n,m}(r).
\end{equation}
The normalization coefficients $N_{n,m}$ are computed to
comply with (\ref{eq.ortho}); the sign of $N_{n,m}$ is chosen
to keep $Q_{n,m}(1)>0$---for compatibility with the Zernike radial functions.

There are some more design choices for compatibility with the Zernike radial functions:
\begin{enumerate}
\item
$Q_{n,n}(r)$ rises monotonously from $r=0$ to $r=1$ with a single local maximum
at the rim $r=1$, 
\item
$Q_{n,m-2}(r)$ has one more extremum than $Q_{n,m}(r)$ in the
range $0\le r\le 1$, 
\item
the lowest-order Taylor coefficient near $r\to 0$ is of degree
\begin{equation}
Q_{n,m}(r)\sim r^{m}.
\label{eq.lowr}
\end{equation}
\end{enumerate}

The final fixture of the $Q$-functions is to 
let the phase functions $\theta_{n,m}(r)$ be polynomials of $r$,
\begin{equation}
\theta_{n,m}(r) = \sum_{i=m,m+2,\ldots n} \beta_{n,m,i}r^i.
\end{equation}
such that
the integrals (\ref{eq.ortho}) are actually 
generalized Fresnel integrals  of chirp functions which are
sums or differences of the radial polynomials $\theta(r)$
\cite{MatharArxiv1211}.

\section{Radial Phase Polynomials}
\subsection{Monotonicity}
The desired number of extrema is established by driving the phase of the sine-functions 
monotonously from 0 up to  multiples of $\pi/2$ as $r$ runs from 0 to 1:
\begin{equation}
\theta_{n,n-2l}(1) = (1+2l)\frac{\pi}{2};\quad m = n-2l>0
\label{eq.r1}
\end{equation}
\begin{equation}
\theta_{n,0}(1) = n\frac{\pi}{2};\quad m=0, n=0,2,4,\ldots
\end{equation}
This clamps the radial functions at the right end-point of the interval.
We could express the polynomials also in an equivalent Bernstein polynomial basis (Appendix \ref{app.Bernst}):
\begin{equation}
\theta_{n,m}(r) \equiv r^m \theta_{n,m}(1)\sum_{i=0}^{(n-m)/2} \alpha_{n,m,i} B_{i,(n-m)/2}(r^2).
\label{eq.alphaDef}
\end{equation}
Dropping a factor $x^m$ and the constant $\theta_{n,m}(1)$ reduces the information
about the radial functions to the expansion coefficients $\alpha$ of ``reduced''
radial phase polynomials $\bar \theta$:
\begin{equation}
\bar \theta_{n,m}(r) \equiv \sum_{i=0}^{(n-m)/2} \alpha_{n,m,i} B_{i,(n-m)/2}(r^2).
\end{equation}

The aim is to find the expansion coefficients $\beta_{n,m,i}$ or $\alpha_{n,m,i}$.
The procedure is to define $Q_{m,m}$ first by finding $\alpha_{m,m,i}$
and then to calculate recursively the $Q_{m+2l,m}$ for $l=1,2,\ldots$ by fixing their
$l$ coefficients $\alpha_{m+2l,m,i}$ by enforcing the orthogonality to the
already established values at smaller $n$.
The $1+(n-m)/2$ degrees of freedom embodied in the coefficients
of (\ref{eq.alphaDef}) match the one degree to clamp the functions at $r=1$
plus $(n-m)/2$ degrees to shape the functions to stay orthogonal with the
already established functions of the same $m$ but smaller $n$.

\subsection{Special Cases}

The piston term at $n=m=0$ has the same constant value as for the Zernike basis:
\begin{equation}
\theta_{0,0}(r)=0 \therefore \cos\theta_{0,0}(r)=1\therefore N_{0,0}=\surd 2
\therefore Q_{0,0}(r) =\surd 2.
\end{equation}

The orthogonality
\begin{equation}
\int_0^1 r \cos(n\frac{\pi}{2}r^2)\cos(n'\frac{\pi}{2}r^2) dr
=\frac12 \int_0^1 \cos(n\frac{\pi}{2}r)\cos(n'\frac{\pi}{2}r) dr
=\frac14 \delta_{n,n'},\quad n-n'\, \textrm{even},\, (n,n')\neq (0,0)
\end{equation}
yields simple phase polynomial coefficients of $m=0$ and all (even) $n$:
\begin{equation}
\beta_{n,0,i} = \left\{ 
\begin{array}{ll}
n\frac{\pi}{2}, &  i=2;\\
0, &  i \neq 2\\
\end{array}
\right.
\label{eq.tn0i}
\end{equation}
\begin{equation}
N_{n,0} = (-1)^{n/2}\times 2.
\label{eq.Nn0}
\end{equation}

Likewise
\begin{equation}
\int_0^1 r \sin(n\frac{\pi}{2}r^2)\sin(n'\frac{\pi}{2}r^2) dr
=\frac12 \int_0^1 \sin(n\frac{\pi}{2}r)\sin(n'\frac{\pi}{2}r) dr
=\frac14 \delta_{n,n'},\quad n-n'\, \textrm{even},\, (n,n')\neq (0,0),
\end{equation}
yields simple phase polynomial coefficients for $m=2$ and all (even) $n$:
\begin{equation}
\beta_{n,2,i} = \left\{ 
\begin{array}{ll}
(n-1)\frac{\pi}{2}, &  i=2;\\
0, &  i \neq 2\\
\end{array}
\right.
\label{eq.tn2i}
\end{equation}
\begin{equation}
N_{n,2} = (-1)^{1+n/2}\times 2.
\label{eq.Nn2}
\end{equation}

For $n=m$ the design choice $\theta \propto r^m$ (\ref{eq.lowr}) and the limit 
(\ref{eq.r1}) on the circle rim enforces
\begin{equation}
\theta_{n,n}(r) = \frac{\pi}{2} r^n = \frac{\pi}{2} B_{n,n}(r)
\therefore
\beta_{n,n,i}=\left\{ \begin{array}{ll} 
\frac{\pi}{2} ,& i=n;\\
0 ,& i\neq n.\\
\end{array}\right.
\label{eq.thetann}
\end{equation}
The normalization constant $N_{n,n}$ allows a semi-analytical 
treatment:
\begin{equation}
\int_0^1 r dr N_{n,n}^2\sin^2(\frac{\pi}{2}r^n)=1.
\label{eq.Nnn}
\end{equation}
A power series is \cite[1.412.1]{GR}
\begin{equation}
\sin^2 x = \frac12[1-\cos(2x)]
=\sum_{k\ge 1}(-1)^{k+1}\frac{2^{2k-1}x^{2k}}{(2k)!}
= x^2-\frac{1}{3}x^4+\frac{2}{45}x^6 -\frac{1}{315}x^8+
\frac{2}{14175}x^{10}-\frac{2}{467775}x^{12}
+\cdots
.
\end{equation}
Insertion into the integrand gives
\begin{equation}
r\sin^2 \frac{\pi r^n}{2} 
=\sum_{k\ge 1}(-1)^{k+1}\frac{\pi^{2k}r^{2kn+1}}{2(2k)!}
= \frac{\pi^2}{2^2}r^{2n+1}-\frac{1}{3}\frac{\pi^4}{2^4}r^{4n+1}+\frac{2}{45}\frac{\pi^6}{2^6}r^{6n+1} -\frac{1}{315}\frac{\pi^8}{2^8}r^{8n+1}+
\frac{2}{14175}\frac{\pi^{10}}{2^{10}}r^{10n+1}
-\cdots .
\end{equation}
According to (\ref{eq.Nnn}), $N_{n,n}$ is the inverse of the square root of the following hypergeometric function:
\begin{eqnarray}
\int_0^1 r\sin^2 \frac{\pi r^n}{2} dr
&=&\sum_{k\ge 1}(-1)^{k+1}\frac{\pi^{2k}}{2(2kn+2)(2k)!}
=\frac{1}{4}-\frac{1}{4}{}_1F_2\left(\frac{1}{n}; 1+\frac{1}{n}, \frac{1}{2}\mid
-\frac{\pi^2}{4}\right)
\nonumber
\\
&=&
\frac{\pi^2}{2^2(2n+2)}
-\frac{1}{3}\frac{\pi^4}{2^4(4n+2)}
+\frac{2}{45}\frac{\pi^6}{2^6(6n+2)}
-\frac{1}{315}\frac{\pi^8}{2^8(8n+2)}
+\frac{2}{14175}\frac{\pi^{10}}{2^{10}(10n+2)}
-\cdots
\end{eqnarray}
In particular
\begin{eqnarray}
N_{1,1} &= & \frac{2\pi}{\sqrt{\pi^2+4}} \approx 1.68712721613; \\
N_{2,2} &= & 2; \\
N_{3,3} &\approx & 2.27799236632; \\
N_{4,4} &\approx & 2.52776603703; \\
N_{5,5} &\approx & 2.75587198375,
\end{eqnarray}
which are useful indicators for error bars in the numerical results
of Appendix \ref{sec.num}.

\subsection{Numerical Synthesis}

The numerical analysis fixes $Q_{n,n}(r)$ via (\ref{eq.thetann})
and rewrites in a loop over $m=n-2, n-4, \ldots$
the orthogonality requirements (\ref{eq.ortho}) and the limiting
value (\ref{eq.r1}) by the
following system of equations (for $m>0$):
\begin{eqnarray}
\int_0^1 r dr \sin\left[\sum_{i=m,m+2,\ldots n} \beta_{n,m,i}r^i\right] \sin\left[\sum_{i=m,m+2,\ldots n-2}\beta_{n-2,m,i}r^i\right] &= & 0 ; \label{eq.betastart}\\
\int_0^1 r dr \sin\left[\sum_{i=m,m+2,\ldots n} \beta_{n,m,i}r^i\right] \sin\left[\sum_{i=m,m+2,\ldots n-4}\beta_{n-4,m,i}r^i\right] &= & 0 ;\\
\ldots &=& 0 ;\\
\int_0^1 r dr \sin\left[\sum_{i=m,m+2,\ldots n} \beta_{n,m,i}r^i\right] \sin\left[\sum_{i=m,m+2,\ldots n-4}\beta_{m,m,i}r^i\right] &= & 0 ;\\
\sum_{i=m,m+2,\ldots n}\beta_{n,m,i} &=& (1+n-m)\frac{\pi}{2} . \label{eq.betaend}
\end{eqnarray}
The $1+(n-m)/2$ unknown values are $\beta_{n,m,m},\beta_{n,m,m+2},\ldots,\beta_{n,m,n}$.
The last equation is used to eliminate $\beta_{n,m,n}$ in all earlier ones and we end up
with a system of $(n-m)/2$ nonlinear equations
\begin{eqnarray}
\int_0^1 r dr \sin\left[\frac{(1+n-m)\pi}{2}r^n+\sum_{i=m,m+2,\ldots n-2} \beta_{n,m,i}(r^i-r^n)
\right] \sin\left[\sum_{i=m,m+2,\ldots n-2}\beta_{n-2,m,i}r^i
\right] &= & 0 ;\\
\int_0^1 r dr \sin\left[\frac{(1+n-m)\pi}{2}r^n+\sum_{i=m,m+2,\ldots n-2} \beta_{n,m,i}(r^i-r^n)
\right] \sin\left[\sum_{i=m,m+2,\ldots n-4}\beta_{n-4,m,i}r^i\right] &= & 0 ;\\
\ldots &=& 0 ;\\
\int_0^1 r dr \sin\left[\frac{(1+n-m)\pi}{2}r^n+\sum_{i=m,m+2,\ldots n-2} \beta_{n,m,i}(r^i-r^n)
\right] \sin\left[\sum_{i=m,m+2,\ldots n-4}\beta_{m,m,i}r^i\right] &= & 0 .
\end{eqnarray}
This system is solved recursively by the common multivariate Newton-method:
insert a set of initial guesses $\{\beta_{n,m,i}\}$,
evaluate the deviations of the right hand side from zero, compute the matrix
of derivatives of the left hand side (the Jacobian)
with respect to
the unknowns (which essentially means replacing $r\to r(r^i-r^n)$ and replacing the first
sine by a cosine), and solve the associated
linear system of equations to calculate updates of the unknown parameters.
All the integrals are evaluated numerically.

The role of the Bernstein basis for the radial phase polynomials is that
the $\alpha$-coefficients are smaller than
the $\beta$-coefficients and essentially of equal sign,
so cancellation of digits is a lesser problem. Instead of (\ref{eq.betastart})-(\ref{eq.betaend}) we may solve
\begin{eqnarray}
\int_0^1 r dr \sin\left[r^m \theta_{n,m}(1) \sum_{i=0}^{\frac{n-m}{2}} \alpha_{n,m,i}B_{i,(n-m)/2}(r^2)\right] 
\sin\left[r^m \theta_{n-2,m}(1) \sum_{i=0}^{\frac{n-m-2}{2}} \alpha_{n-2,m,i}B_{i,(n-m-2)/2}(r^2)\right] &= & 0 ;\\
\int_0^1 r dr \sin\left[r^m \theta_{n,m}(1) \sum_{i=0}^{\frac{n-m}{2}} \alpha_{n,m,i}B_{i,(n-m)/2}(r^2)\right] 
\sin\left[r^m \theta_{n-4,m}(1) \sum_{i=0}^{\frac{n-m-4}{2}} \alpha_{n-4,m,i}B_{i,(n-m-4)/2}(r^2)\right] &= & 0 ;\\
\ldots &=& 0 ;\\
\int_0^1 r dr \sin\left[r^m \theta_{n,m}(1) \sum_{i=0}^{\frac{n-m}{2}} \alpha_{n,m,i}B_{i,(n-m)/2}(r^2)\right] 
\sin\left[r^m \theta_{m,m}(1) \sum_{i=0}^{0} \alpha_{m,m,i}B_{i,0}(r^2)\right] &= & 0 ;\\
\alpha_{n,m,(n-m)/2} &=& 1 ;
\end{eqnarray}

The last equation is inserted into the initial set of equations, so the
orthogonality requirements are:
\begin{eqnarray}
\int_0^1 r dr \sin\left[r^n\theta_{n,m}(1)+r^m \theta_{n,m}(1) 
\sum_{i=0}^{\frac{n-m}{2}-1}\!\!\! \alpha_{n,m,i}B_{i,\frac{n-m}{2}}(r^2)\right]
\sin\left[r^m \theta_{n-2,m}(1) \sum_{i=0}^{\frac{n-m-2}{2}}\!\!\! 
\alpha_{n-2,m,i}B_{i,\frac{n-m-2}{2}}(r^2)\right] &= & 0 ; \label{eq.alphA}\\
\int_0^1 r dr \sin\left[r^n\theta_{n,m}(1)+r^m \theta_{n,m}(1) 
\sum_{i=0}^{\frac{n-m}{2}-1}\!\! \alpha_{n,m,i}B_{i,\frac{n-m}{2}}(r^2)\right]
\sin\left[r^m \theta_{n-4,m}(1) \sum_{i=0}^{\frac{n-m-4}{2}}\!\!\! 
\alpha_{n-4,m,i}B_{i,\frac{n-m-4}{2}}(r^2)\right] &= & 0 ;\\
\ldots &=& 0 ;\\
\int_0^1 r dr \sin\left[r^n\theta_{n,m}(1)+r^m \theta_{n,m}(1) 
\sum_{i=0}^{\frac{n-m}{2}-1}\!\!\! \alpha_{n,m,i}B_{i,\frac{n-m}{2}}(r^2)\right]
\sin\left[r^m \theta_{m,m}(1) \sum_{i=0}^{0} \alpha_{m,m,i}B_{i,0}(r^2)\right] &= & 0 . \label{eq.alphE}
\end{eqnarray}
The unknown values $\alpha_{n,m,0},\alpha_{n,m,1},\ldots,\alpha_{n,m,(n-m-2)/2}$ are extracted
with again with the multivariate Newton method.

\subsection{Initial values}
A numerical problem here is that for most guesses for the $\beta_{n,m,i}$
this Newton procedure branches into parameter
sets that yield non-monotonous $Q_{n,m}(r)$, which fulfill the orthogonality relations
but define radial functions that violate the minimax principle. 
Two rules of thumb to start with guesses in the attractor region are:
\begin{enumerate}
\item
A good initial guess proposed by heuristics of the expansion coefficients of the Bernstein polynomials is
\begin{equation}
\alpha_{n,n-2,0}\approx 3-\pi +(\frac{\pi}{4}-\frac{1}{2})n
\approx -\frac{1}{6}+\frac{7}{24} n,
\end{equation}
so
\begin{equation}
\beta_{n,n-2,n-2}=\frac{3}{2}\pi \alpha_{n,n-2,0},\quad
\beta_{n,n-2,n}=\frac{3}{2}\pi (1-\alpha_{n,n-2,0}).
\end{equation}
\item
For systems with 3 or more unknown values one may ``lift'' an already
known coefficient set of a smaller $n$ as
\begin{eqnarray}
\beta_{n,m,m} &\approx&  \pi+\beta_{n-2,m,m} ; \\
\beta_{n,m,i} &\approx&  \beta_{n-2,m,i} ,\quad m+2\le i \le n-2; \\
\beta_{n,m,n} &\approx&  0 .
\end{eqnarray}
Although these guesses are far off the final, converged values, they
seem to steer towards the correct phase polynomials if sufficient damping
is included in the Newton updates.
\end{enumerate}

\begin{figure}
\includegraphics[width=0.45\columnwidth]{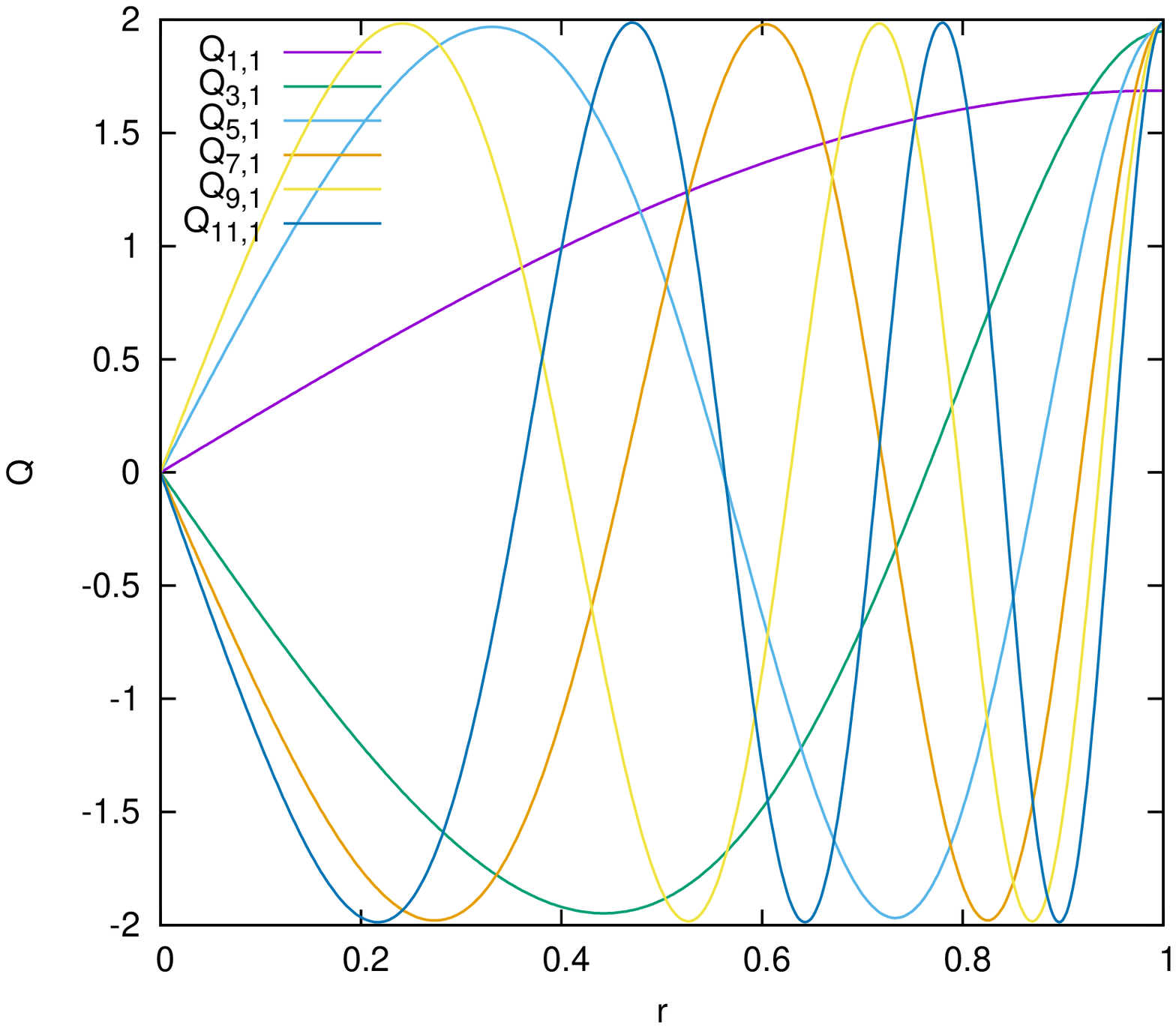}
\includegraphics[width=0.45\columnwidth]{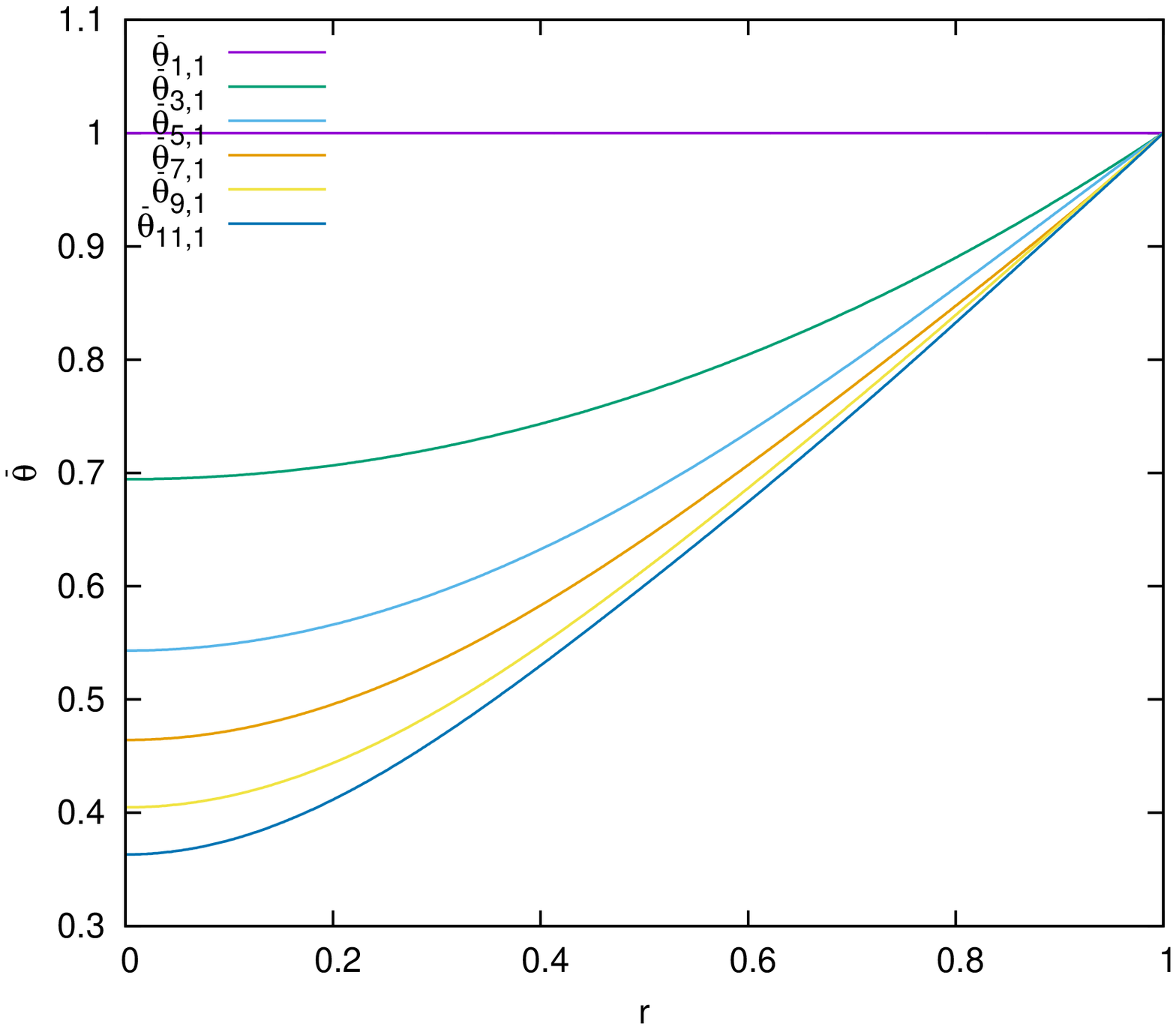}
\caption{Radial functions $Q_{n,1}(r)$ and associated reduced polynomials $\bar \theta_{n,1}(r)$.
The ``tip-tilt'' representative $Q_{1,1}$ is no longer a straight line as it used to be in the Zernike case.
}
\label{fig.Q1}
\end{figure}

\begin{figure}
\includegraphics[width=0.45\columnwidth]{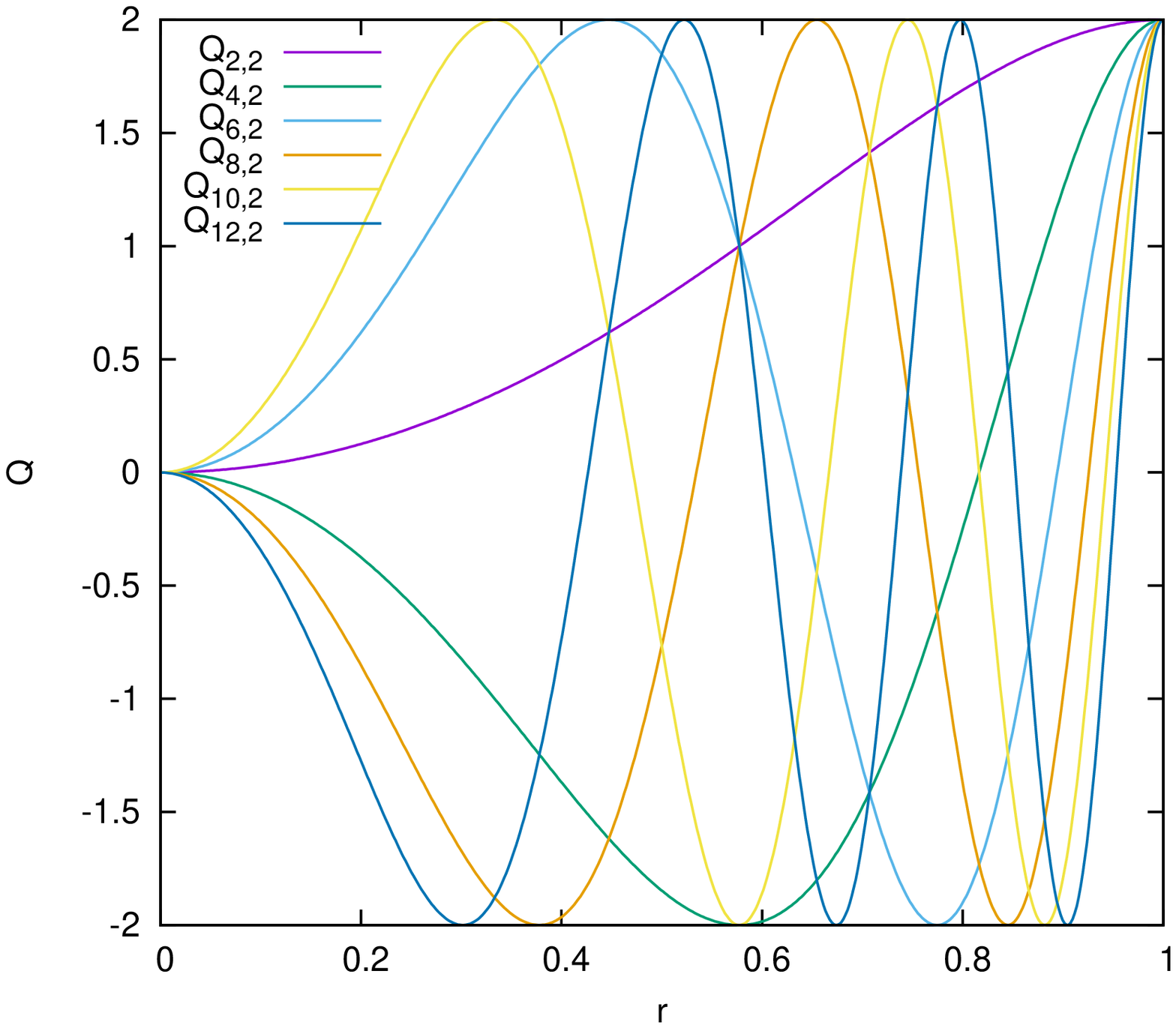}
\includegraphics[width=0.45\columnwidth]{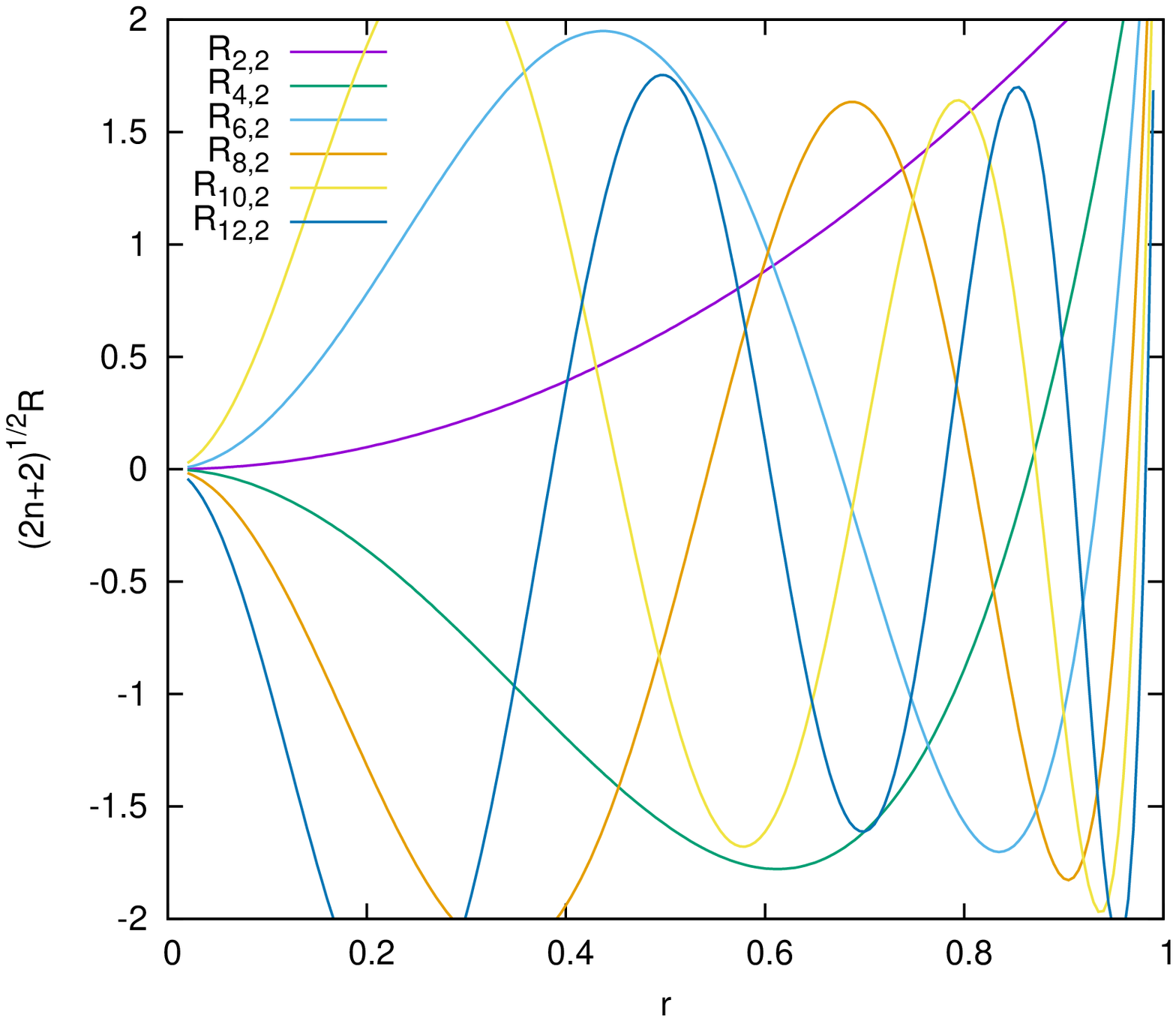}
\caption{Radial functions $Q_{n,2}(r)$ in comparison with normalized Zernike polynomials $\sqrt{2n+2}R_n^{2}(r)$.}
\label{fig.Q2}
\end{figure}

\begin{figure}
\includegraphics[width=0.45\columnwidth]{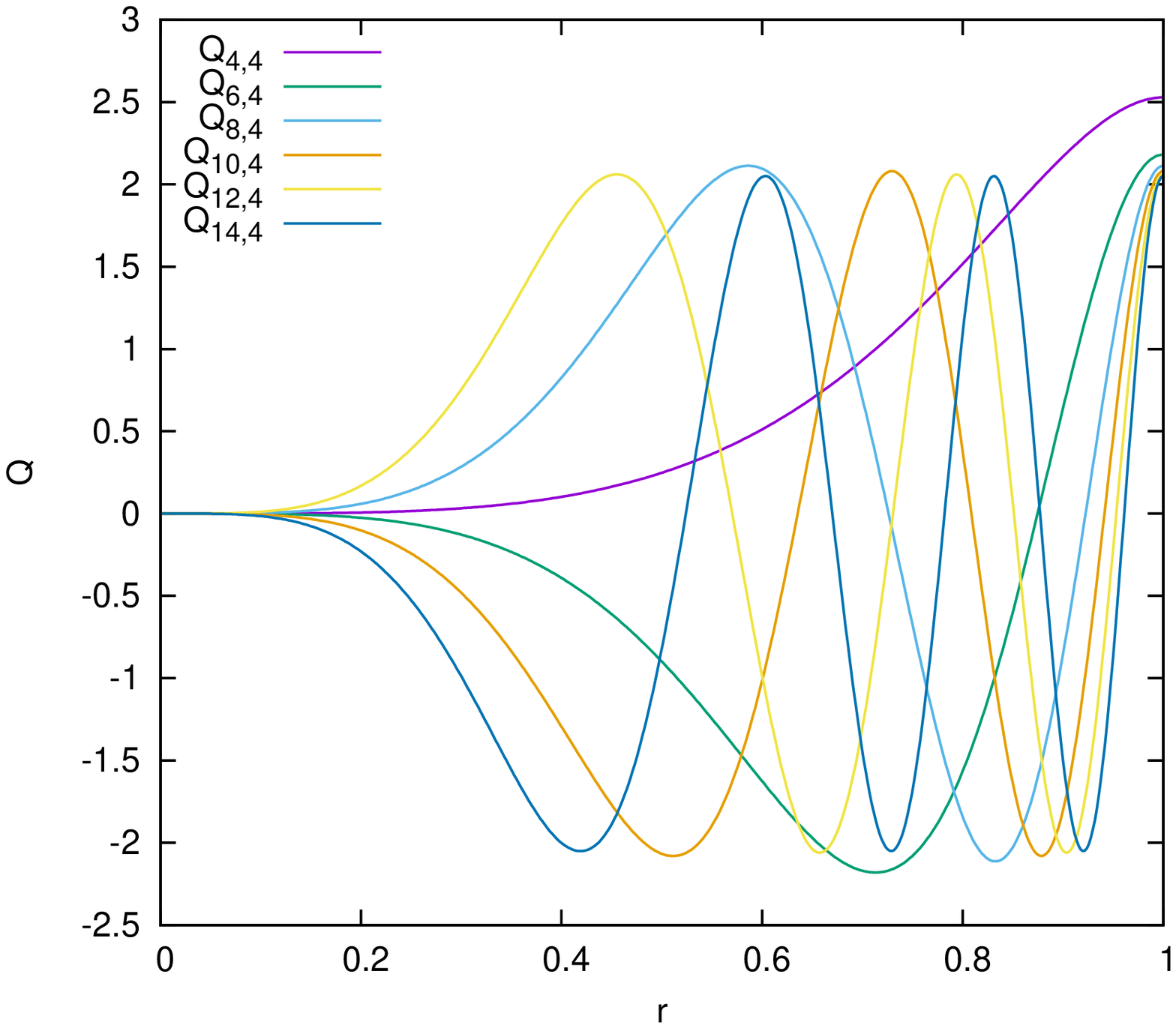}
\includegraphics[width=0.45\columnwidth]{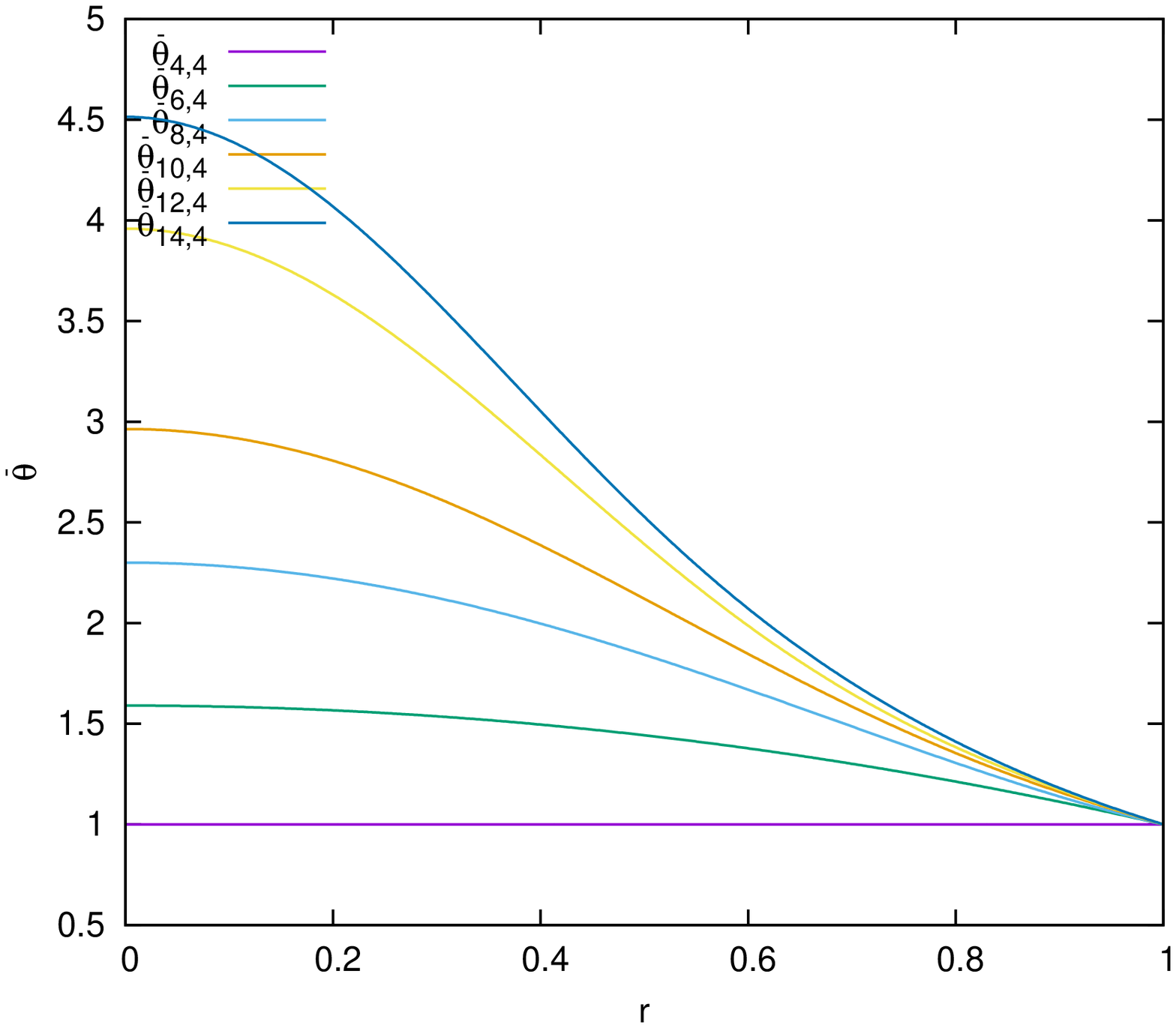}
\caption{Radial functions $Q_{n,4}(r)$ and associated reduced phase polynomials $\bar \theta_{n,4}(r)$.}
\label{fig.Q4}
\end{figure}

\begin{figure}
\includegraphics[width=0.45\columnwidth]{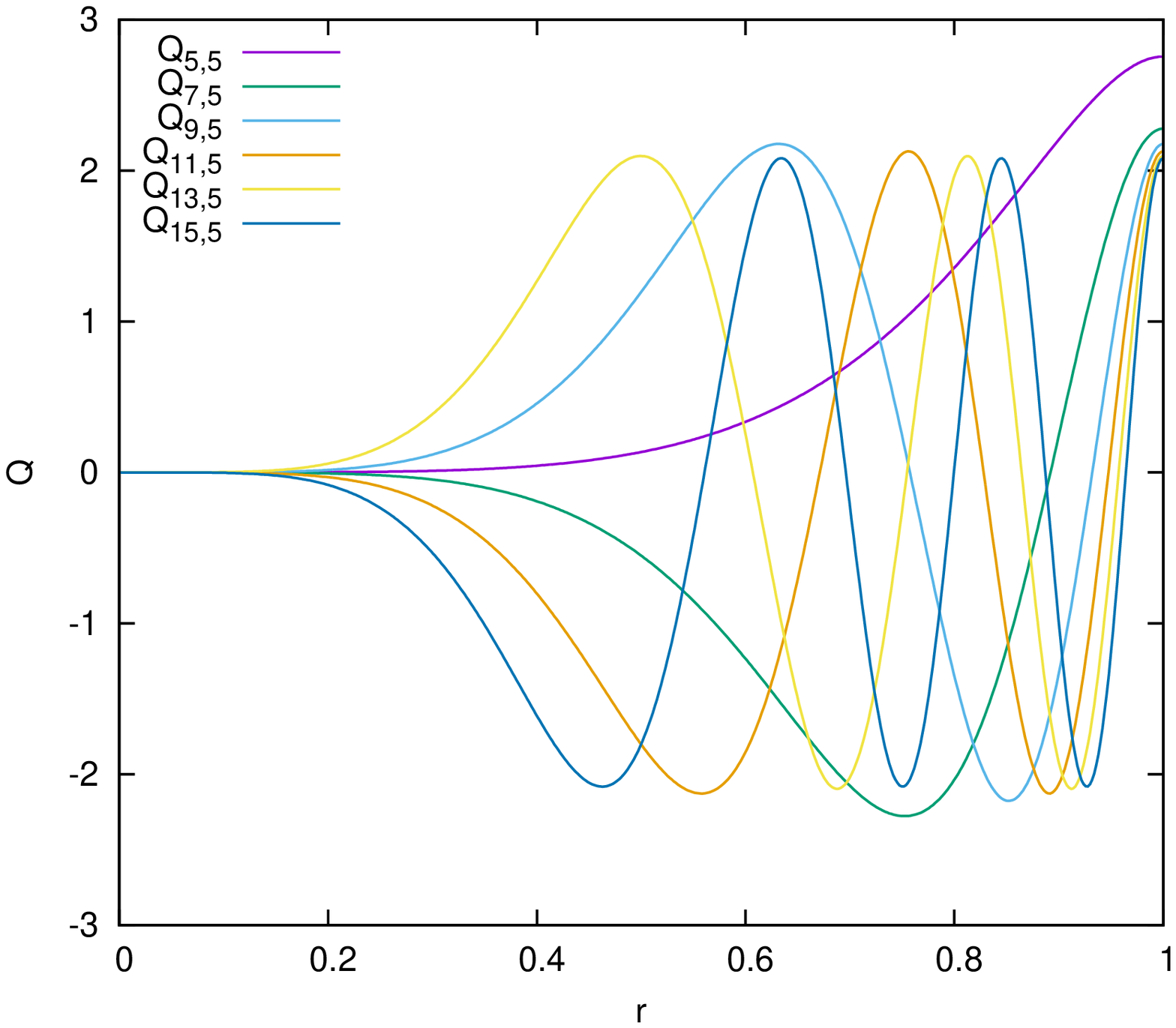}
\includegraphics[width=0.45\columnwidth]{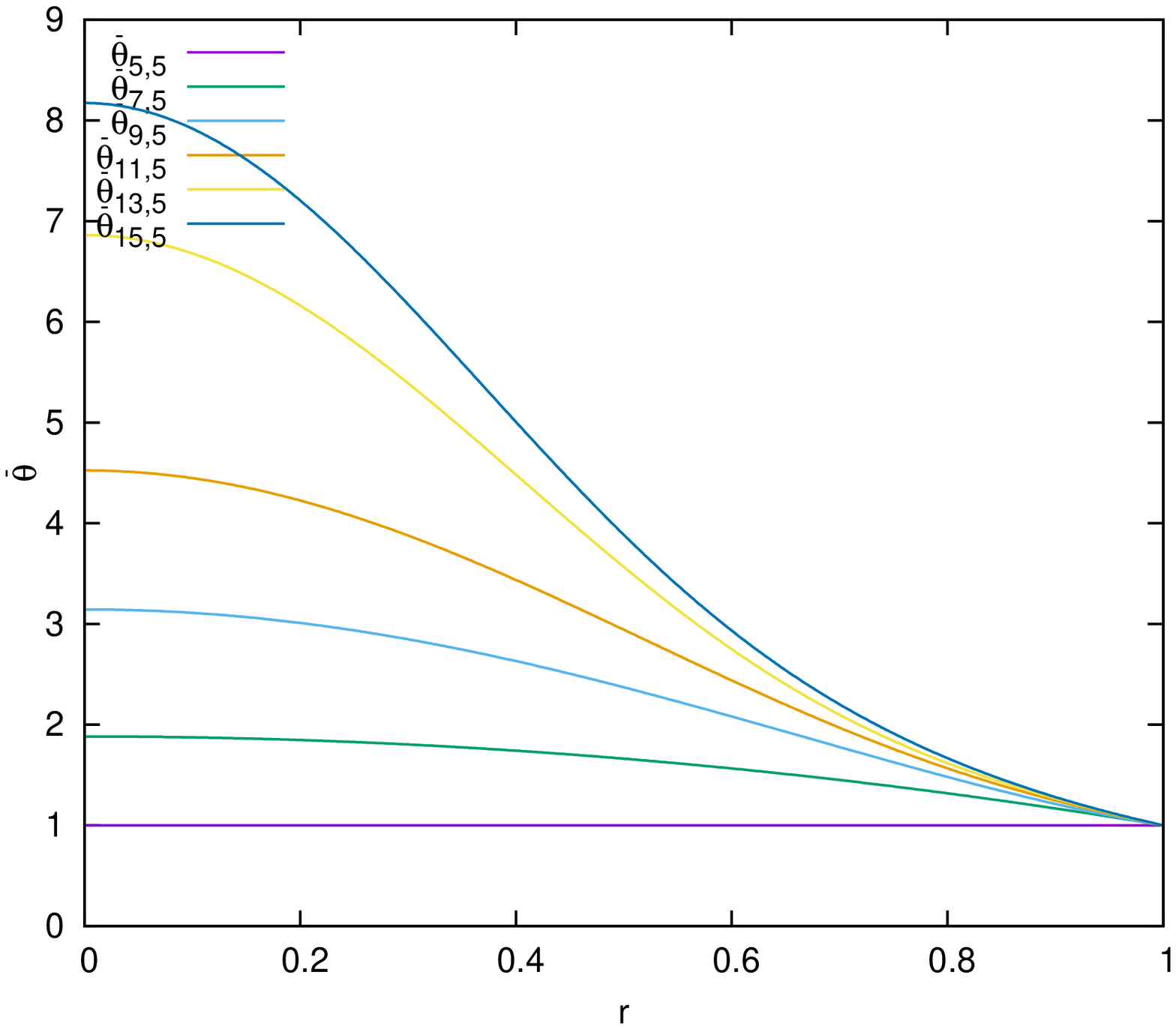}
\caption{Radial functions $Q_{n,5}(r)$ and associated reduced phase polynomials $\bar \theta_{n,5}(r)$.}
\label{fig.Q5}
\end{figure}

The first few radial functions for azimuthal quantum numbers
$m=1,4$ and $5$ are plotted in Figures \ref{fig.Q1}--\ref{fig.Q5}.

\begin{figure}
\includegraphics[width=0.45\columnwidth]{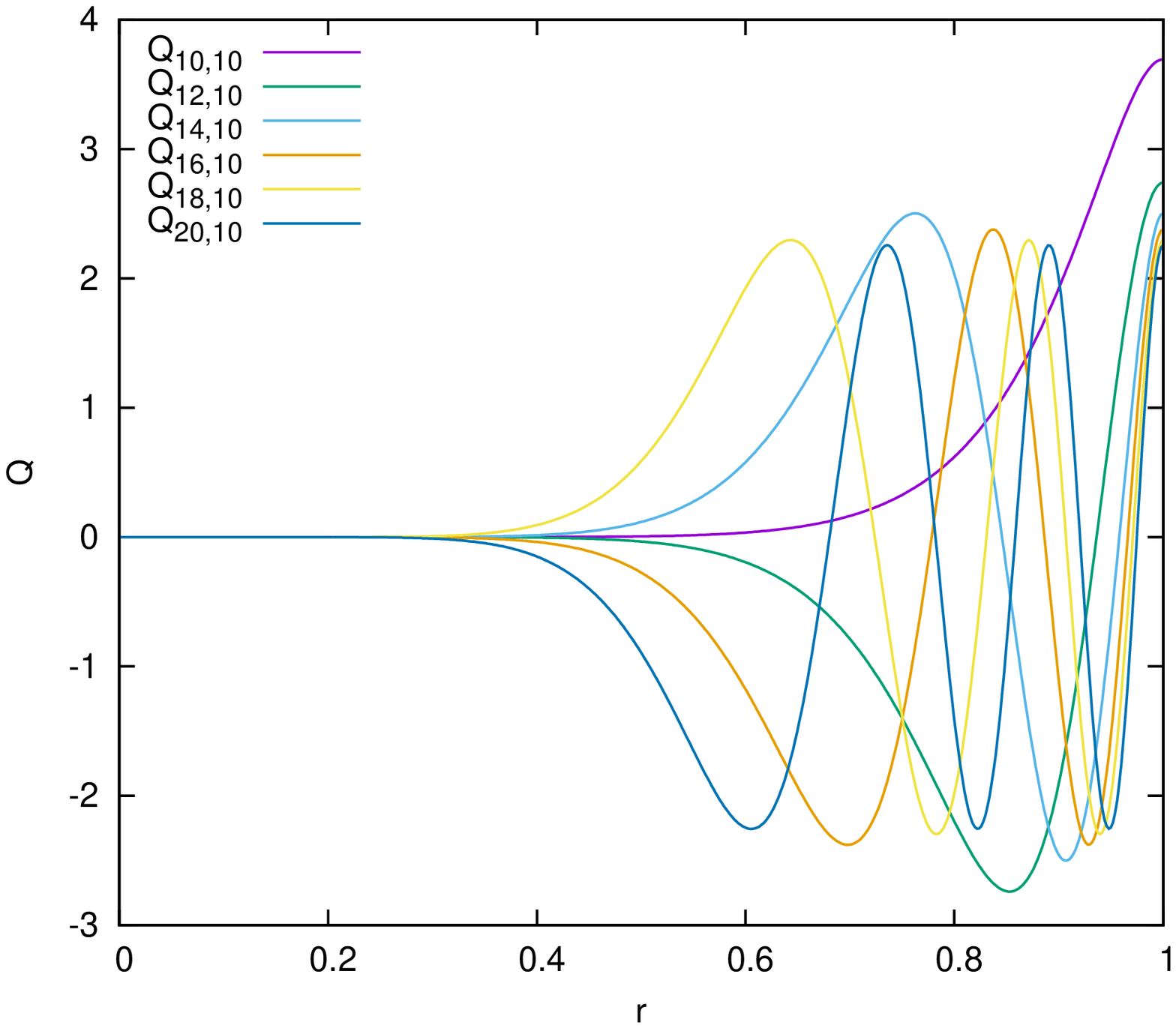}
\includegraphics[width=0.45\columnwidth]{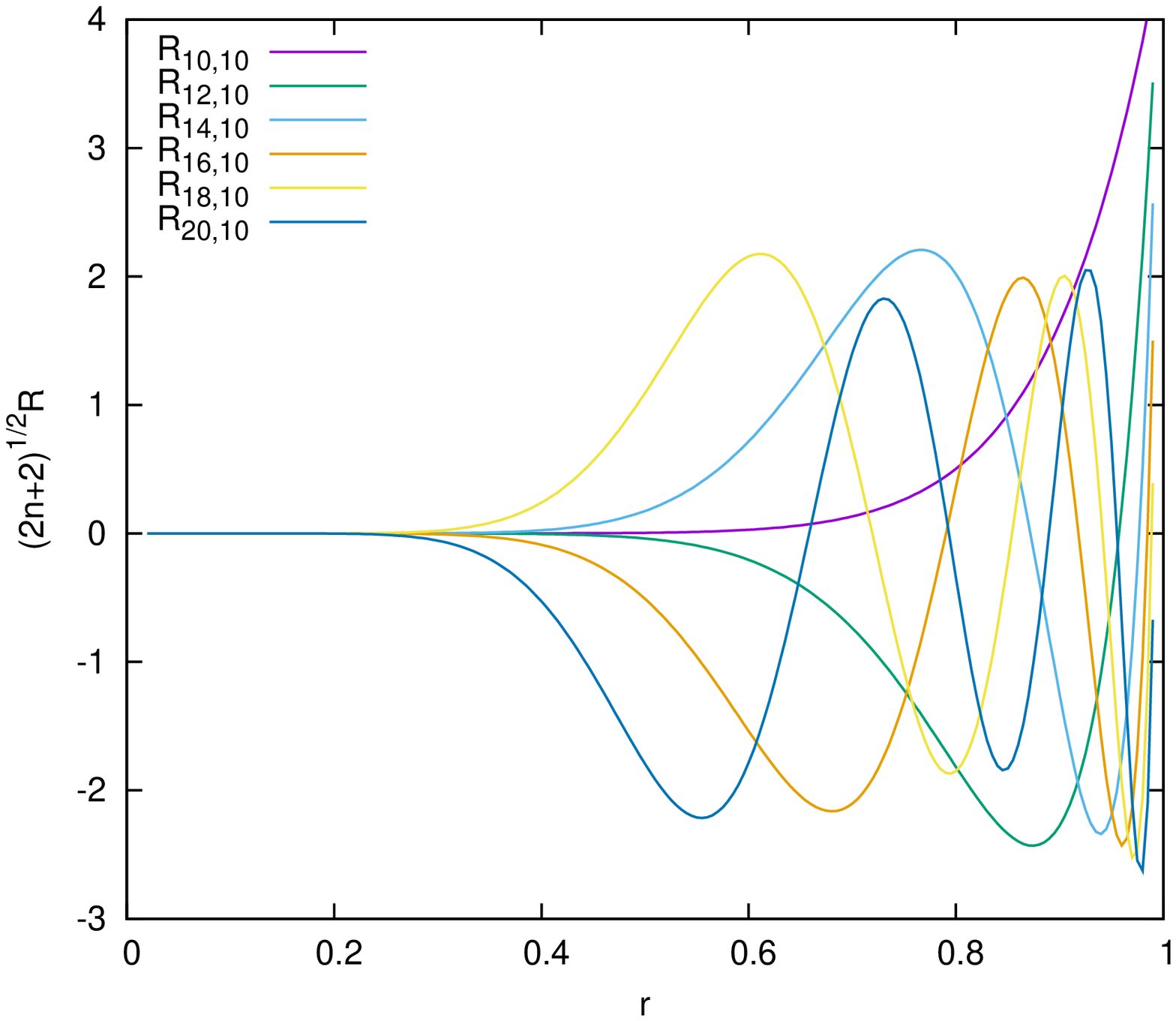}
\caption{Radial functions $Q_{n,10}(r)$ in comparison with normalized Zernike polynomials $\sqrt{2n+2}R_n^{10}(r)$.}
\label{fig.Q10}
\end{figure}

\appendix

\section{Expansion Coefficients}\label{sec.num}

The results of the numerical procedure are summarized complete for $n \le 22$
in the following table for all $0<m\le n$.
The coefficients for $m=0$ are precisely known
by Eqs.\ (\ref{eq.tn0i}) and (\ref{eq.Nn0}) and not listed.
The coefficients for $m=2$ are precisely known
by Eqs.\ (\ref{eq.tn2i}) and (\ref{eq.Nn2}) but are run through
the same numerical procedure to produce indicators of error bars.

The estimated accuracy of the numerical coefficients
is of the order of four digits, although the values
are printed to much larger apparent accuracy.
They have been computed by evaluating the integrals
(\ref{eq.alphA})--(\ref{eq.alphE}) with a Romberg rule over 
4096 sampling points in the interval $0\le r\le 1$.
Superior Gauss-Legendre methods \cite{SwarztrauberSIAMJ24,HaleOccam,HaleSIAMJSci35} have not been employed.

An entry contains three types of lines. 
\begin{enumerate}
\item
A value of $n$, a value of $m$, the letter \texttt{N} and the signed normalization constant $N_{n,m}$;
\item
A value of $n$, a value of $m$, a value of $i$, the letter \texttt{t} and the coefficient $\beta_{n,m,i}$;
\item
A value of $n$, a value of $m$, a value of $i$, the letter \texttt{B} and the coefficient $\alpha_{n,m,i}$.
\end{enumerate}
The fields in the line are separated by a blank.
To save some paper, three columns of these entries have been
printed in a single line of the print.

\small
\begin{Verbatim}
1 1 N 1.6871271565051231             1 1 0 B 1.0000000000000000           1 1 1 t 1.5707963267948966
\end{Verbatim}
\begin{Verbatim}
2 2 N 1.9999999006588837             2 2 0 B 1.0000000000000000           2 2 2 t 1.5707963267948966
\end{Verbatim}
\begin{Verbatim}
3 3 N 2.2779922195323444             3 1 N -1.9479336816837157            3 1 1 t 3.2722064718377821
3 3 0 B 1.0000000000000000           3 1 0 B 0.6943837797470318           3 1 3 t 1.4401825085469078
3 3 3 t 1.5707963267948966           3 1 1 B 1.0000000000000000
\end{Verbatim}
\begin{Verbatim}
4 4 N 2.5277658364697498             4 2 N -1.9999999006584910            4 2 2 t 4.7123889803846899
4 4 0 B 1.0000000000000000           4 2 0 B 1.0000000000000000           4 2 4 t 0.0000000000000000
4 4 4 t 1.5707963267948966           4 2 1 B 1.0000000000000000
\end{Verbatim}
\begin{Verbatim}
5 5 N 2.7558717238509236             5 3 1 B 1.0000000000000000           5 1 1 B 0.8323536660659460
5 5 0 B 1.0000000000000000           5 3 3 t 6.1094671372212254           5 1 2 B 1.0000000000000000
5 5 5 t 1.5707963267948966           5 3 5 t -1.3970781568365356          5 1 1 t 4.2653496600138467
5 3 N -2.0854065890250237            5 1 N 1.9683580365897798             5 1 3 t 4.5438814924788462
5 3 0 B 1.2964691927283318           5 1 0 B 0.5430811859252311           5 1 5 t -0.9552495185182099
\end{Verbatim}
\begin{Verbatim}
6 6 N 2.9669226169055586             6 4 1 B 1.0000000000000000           6 2 1 B 1.0000000000000000
6 6 0 B 1.0000000000000000           6 4 4 t 7.4925051767445798           6 2 2 B 1.0000000000000000
6 6 6 t 1.5707963267948966           6 4 6 t -2.7801161963598900          6 2 2 t 7.8539816339744831
6 4 N -2.1804596649343622            6 2 N 1.9999999006577055             6 2 4 t 0.0000000000000000
6 4 0 B 1.5899589800273531           6 2 0 B 1.0000000000000000           6 2 6 t 0.0000000000000000
\end{Verbatim}
\begin{Verbatim}
7 7 N 3.1641669053471717             7 3 N 2.0522748340685570             7 1 0 B 0.4642215522790204
7 7 0 B 1.0000000000000000           7 3 0 B 1.5855622286135586           7 1 1 B 0.7326398213825204
7 7 7 t 1.5707963267948966           7 3 1 B 1.1510886537962158           7 1 2 B 0.8620164173719093
7 5 N -2.2777338754943961            7 3 2 B 1.0000000000000000           7 1 3 B 1.0000000000000000
7 5 0 B 1.8821433358763838           7 3 3 t 12.4529766230545399          7 1 1 t 5.1043825639723723
7 5 1 B 1.0000000000000000           7 3 5 t -6.8246949541252974          7 1 3 t 8.8542390542008598
7 5 5 t 8.8693915154883513           7 3 7 t 2.2256999650452406           7 1 5 t -4.5865291373803108
7 5 7 t -4.1570025351036614          7 1 N -1.9796188231171719            7 1 7 t 1.6234818067713550
\end{Verbatim}
\begin{Verbatim}
8 8 N 3.3499625464632930             8 4 N 2.1128329254446929             8 2 0 B 1.0000000000000000
8 8 0 B 1.0000000000000000           8 4 0 B 2.3001652356752047           8 2 1 B 1.0000000000000000
8 8 8 t 1.5707963267948966           8 4 1 B 1.2965456834639862           8 2 2 B 1.0000000000000000
8 6 N -2.3744207306006835            8 4 2 B 1.0000000000000000           8 2 3 B 1.0000000000000000
8 6 0 B 2.1736490254813351           8 4 4 t 18.0654555160996464          8 2 2 t 10.9955742875642763
8 6 1 B 1.0000000000000000           8 4 6 t -15.7648190611292096         8 2 4 t 0.0000000000000000
8 6 6 t 10.2430797149021633          8 4 8 t 5.5533451790040463           8 2 6 t 0.0000000000000000
8 6 8 t -5.5306907345174734          8 2 N -1.9999999006565273            8 2 8 t 0.0000000000000000
\end{Verbatim}
\begin{Verbatim}
9 9 N 3.5260753645649359             9 5 5 t 24.6905274862244463          9 1 N 1.9833803657606531
9 9 0 B 1.0000000000000000           9 5 7 t -26.7691103496335822         9 1 0 B 0.4047043188429115
9 9 9 t 1.5707963267948966           9 5 9 t 9.9325644973836191           9 1 1 B 0.6577174198136432
9 7 N -2.4693716976346138            9 3 N -2.0364258316568831            9 1 2 B 0.7719278551280649
9 7 0 B 2.4647599453238663           9 3 0 B 1.8113117427617004           9 1 3 B 0.9015499603695294
9 7 1 B 1.0000000000000000           9 3 1 B 1.2899888786536680           9 1 4 B 1.0000000000000000
9 7 7 t 11.6149076056377582          9 3 2 B 1.1281425361877348           9 1 1 t 5.7213725172882839
9 7 9 t -6.9025186252530683          9 3 3 B 1.0000000000000000           9 1 3 t 14.3075537868892190
9 5 N 2.1770017277838954             9 3 3 t 19.9164128254737921          9 1 5 t -11.7736587375626823
9 5 0 B 3.1436955975831435           9 3 5 t -17.1967328403169417         9 1 7 t 8.7206152274440650
9 5 1 B 1.4395211038565038           9 3 7 t 11.8579523950527310          9 1 9 t -2.8387158529048160
9 5 2 B 1.0000000000000000           9 3 9 t -3.5820580926453051
\end{Verbatim}
\begin{Verbatim}
10 10 N 3.6938645796605060           10 6 6 t 32.3273155698162636         10 2 N 1.9999999006549564
10 10 0 B 1.0000000000000000         10 6 8 t -39.8169994932505851        10 2 0 B 1.0000000000000000
10 10 10 t 1.5707963267948966        10 6 10 t 15.3436655574088046        10 2 1 B 1.0000000000000000
10 8 N -2.5621266646518794           10 4 N -2.0803386368920199           10 2 2 B 1.0000000000000000
10 8 0 B 2.7556222038343123          10 4 0 B 2.9638587431473974          10 2 3 B 1.0000000000000000
10 8 1 B 1.0000000000000000          10 4 1 B 1.6087736717897403          10 2 4 B 1.0000000000000000
10 8 8 t 12.9855637074521867         10 4 2 B 1.2507914828616563          10 2 2 t 14.1371669411540696
10 8 10 t -8.2731747270674969        10 4 3 B 1.0000000000000000          10 2 4 t 0.0000000000000000
10 6 N 2.2425818193058856            10 4 4 t 32.5893289881240961         10 2 6 t 0.0000000000000000
10 6 0 B 4.1160416558622796          10 4 6 t -44.6998157042473730        10 2 8 t 0.0000000000000000
10 6 1 B 1.5812127404869507          10 4 8 t 32.8911564482965201         10 2 10 t 0.0000000000000000
10 6 2 B 1.0000000000000000          10 4 10 t -9.7850954446089668
\end{Verbatim}
\begin{Verbatim}
11 11 N 3.8544012493346635           11 5 0 B 4.5245441521751427          11 3 7 t 49.2300531493700388
11 11 0 B 1.0000000000000000         11 5 1 B 1.9561144488077112          11 3 9 t -33.5086177004263925
11 11 11 t 1.5707963267948966        11 5 2 B 1.3704880380148920          11 3 11 t 9.2974816710964301
11 9 N -2.6525363593996709           11 5 3 B 1.0000000000000000          11 1 N -1.9873005399434828
11 9 0 B 3.0463182443855785          11 5 5 t 49.7499613426063069         11 1 0 B 0.3630829521276220
11 9 1 B 1.0000000000000000          11 5 7 t -84.7240788172898108        11 1 1 B 0.6195795042127378
11 9 9 t 14.3554365255874349         11 5 9 t 65.4061827033935778         11 1 2 B 0.6902653057534245
11 9 11 t -9.6430475452027451        11 5 11 t -19.4364909411457975       11 1 3 B 0.8303873542803619
11 7 N 2.3084287008014712            11 3 N 2.0275990551056175            11 1 4 B 0.9128914926150938
11 7 0 B 5.2171257824976542          11 3 0 B 2.1011679756827771          11 1 5 B 1.0000000000000000
11 7 1 B 1.7221641428981020          11 3 1 B 1.3834443894886579          11 1 1 t 6.2736230427630761
11 7 2 B 1.0000000000000000          11 3 2 B 1.2461064556944625          11 1 3 t 22.1597113017970598
11 7 7 t 40.9752100778713298         11 3 3 B 1.0965917139446977          11 1 5 t -32.1057928877611090
11 7 9 t -54.8987290577204591        11 3 4 B 1.0000000000000000          11 1 7 t 44.1035150761279513
11 7 11 t 21.7775006138236125        11 3 3 t 29.7045624436341749         11 1 9 t -33.0284487250607014
11 5 N -2.1279685890228062           11 3 5 t -40.5863126225201818        11 1 11 t 9.8761517868775860
\end{Verbatim}
\begin{Verbatim}
12 12 N 4.0085463540699856           12 6 0 B 6.5599173686526489          12 4 8 t 159.7982867342640598
12 12 0 B 1.0000000000000000         12 6 1 B 2.3314556048485304          12 4 10 t -109.5482315108087225
12 12 12 t 1.5707963267948966        12 6 2 B 1.4885201221016779          12 4 12 t 30.1016013063744597
12 10 N -2.7405982724390614          12 6 3 B 1.0000000000000000          12 2 N -1.9999999006529925
12 10 0 B 3.3368979406950816         12 6 6 t 72.1300587473033723         12 2 0 B 1.0000000000000000
12 10 1 B 1.0000000000000000         12 6 8 t -139.4830963380997620       12 2 1 B 1.0000000000000000
12 10 10 t 15.7247610843998668       12 6 10 t 111.6774171775991475       12 2 2 B 1.0000000000000000
12 10 12 t -11.0123721040151769      12 6 12 t -33.3288052992384815       12 2 3 B 1.0000000000000000
12 8 N 2.3739159231709517            12 4 N 2.0606352607461742            12 2 4 B 1.0000000000000000
12 8 0 B 6.4468959183547374          12 4 0 B 3.9588705312318715          12 2 5 B 1.0000000000000000
12 8 1 B 1.8626539739617522          12 4 1 B 1.7982236403856168          12 2 2 t 17.2787595947438628
12 8 2 B 1.0000000000000000          12 4 2 B 1.5214795074189302          12 2 4 t 0.0000000000000000
12 8 8 t 50.6338021389031660         12 4 3 B 1.1914000710073076          12 2 6 t 0.0000000000000000
12 8 10 t -72.0091040739159577       12 4 4 B 1.0000000000000000          12 2 8 t 0.0000000000000000
12 8 12 t 29.2292835689872748        12 4 4 t 55.9672135984402638         12 2 10 t 0.0000000000000000
12 6 N -2.1774566893976795           12 4 6 t -122.1817031871159913       12 2 12 t 0.0000000000000000
\end{Verbatim}
\begin{Verbatim}
13 13 N 4.1570038942323646           13 7 9 t -211.1676232605324680       13 3 3 t 39.2228367588695030
13 13 0 B 1.0000000000000000         13 7 11 t 173.9287186915836537       13 3 5 t -67.9623471371769554
13 13 13 t 1.5707963267948966        13 7 13 t -52.2217248299338382       13 3 7 t 117.2719326708493161
13 11 N -2.8263796837859086          13 5 N 2.0969017565107164            13 3 9 t -132.8055037768876705
13 11 0 B 3.6273931556552369         13 5 0 B 6.8624432428924186          13 3 11 t 83.0123864710731159
13 11 1 B 1.0000000000000000         13 5 1 B 2.2391778606081654          13 3 13 t -21.4605453919834464
13 11 11 t 17.0936875342325842       13 5 2 B 1.8287856347248273          13 1 N 1.9887853985258073
13 11 13 t -12.3812985538478943      13 5 3 B 1.2845053729414280          13 1 0 B 0.3288884690286226
13 9 N 2.4386936298203388            13 5 4 B 1.0000000000000000          13 1 1 B 0.5748036169612312
13 9 0 B 7.8053166528499373          13 5 5 t 97.0155057489648265         13 1 2 B 0.6293977595147509
13 9 1 B 2.0028384608675147          13 5 7 t -261.4394980904439000       13 1 3 B 0.7856396792121538
13 9 2 B 1.0000000000000000          13 5 9 t 357.3485466836784794        13 1 4 B 0.8403494650788614
13 9 9 t 61.3028136388385941         13 5 11 t -245.8035545155552256      13 1 5 B 0.9292045920186192
13 9 11 t -91.1451143027348246       13 5 13 t 67.0161671145098893        13 1 6 B 1.0000000000000000
13 9 13 t 37.6962822978707136        13 3 N -2.0225187135941081           13 1 1 t 6.7160183879796467
13 7 N -2.2277977712971895           13 3 0 B 2.2700030371856641          13 1 3 t 30.1300436639096905
13 7 0 B 9.1360579319681757          13 3 1 B 1.4833453287486434          13 1 5 t -58.6026348351950150
13 7 1 B 2.7344634438003968          13 3 2 B 1.3753933574209175          13 1 7 t 119.6505147376909828
13 7 2 B 1.6055575890431455          13 3 3 B 1.1775415815333572          13 1 9 t -151.9729663315577781
13 7 3 B 1.0000000000000000          13 3 4 B 1.0820449426696382          13 1 11 t 102.3067096221325469
13 7 7 t 100.4562036864469288        13 3 5 B 1.0000000000000000          13 1 13 t -27.8073329966264178
\end{Verbatim}
\begin{Verbatim}
14 14 N 4.3003580200619583           14 8 10 t -301.9517635933222994      14 4 4 t 78.0020428444088330
14 14 0 B 1.0000000000000000         14 8 12 t 254.3517984506937032       14 4 6 t -206.2140645485721882
14 14 14 t 1.5707963267948966        14 8 14 t -76.8557276822191525       14 4 8 t 357.2286848144291758
14 12 N -2.9099796214415765          14 6 N 2.1351776326140948            14 4 10 t -377.7958077270740428
14 12 0 B 3.9178251598979772         14 6 0 B 11.1410985911011365         14 4 12 t 219.1653093487865425
14 12 1 B 1.0000000000000000         14 6 1 B 2.7037817077344996          14 4 14 t -53.1074051372344576
14 12 12 t 18.4623161105771132       14 6 2 B 2.1679815856954638          14 2 N 1.9999999006506359
14 12 14 t -13.7499271301924233      14 6 3 B 1.3762830243400433          14 2 0 B 1.0000000000000000
14 10 N 2.5025680363719200           14 6 4 B 1.0000000000000000          14 2 1 B 1.0000000000000000
14 10 0 B 9.2923631168567868         14 6 6 t 157.5035706902531676        14 2 2 B 1.0000000000000000
14 10 1 B 2.1428113858676431         14 6 8 t -477.1190292622876213       14 2 3 B 1.0000000000000000
14 10 2 B 1.0000000000000000         14 6 10 t 670.2303692592919979       14 2 4 B 1.0000000000000000
14 10 10 t 72.9820492560150873       14 6 12 t -461.2909619992473954      14 2 5 B 1.0000000000000000
14 10 12 t -112.3048959726784183     14 6 14 t 124.8132182531439208       14 2 6 B 1.0000000000000000
14 10 14 t 47.1768283506378141       14 4 N -2.0506554559690726           14 2 2 t 20.4203522483336560
14 8 N -2.2784094620367242           14 4 0 B 4.5143311599830799          14 2 4 t 0.0000000000000000
14 8 0 B 12.3187078337148858         14 4 1 B 2.1274229630393389          14 2 6 t 0.0000000000000000
14 8 1 B 3.1649715001543821          14 4 2 B 1.8079588025477104          14 2 8 t 0.0000000000000000
14 8 2 B 1.7219677394970000          14 4 3 B 1.3694633956296814          14 2 10 t 0.0000000000000000
14 8 3 B 1.0000000000000000          14 4 4 B 1.1622790839248309          14 2 12 t 0.0000000000000000
14 8 8 t 135.4512671124120250        14 4 5 B 1.0000000000000000          14 2 14 t 0.0000000000000000
\end{Verbatim}
\begin{Verbatim}
15 15 N 4.4390996495054235           15 7 2 B 2.5384440735242009          15 3 5 B 1.0687546033272857
15 15 0 B 1.0000000000000000         15 7 3 B 1.4670932446284237          15 3 6 B 1.0000000000000000
15 15 15 t 1.5707963267948966        15 7 4 B 1.0000000000000000          15 3 3 t 51.6579454603501970
15 13 N -2.9915094406436906          15 7 7 t 242.6661063167431191        15 3 5 t -118.4515848858738342
15 13 0 B 4.2082086756364336         15 7 9 t -790.2534416611649117       15 3 7 t 271.1996533582855461
15 13 1 B 1.0000000000000000         15 7 11 t 1130.0821329122048439      15 3 9 t -419.5780232946832802
15 13 13 t 19.8307161902283793       15 7 13 t -777.1061976376817831      15 3 11 t 395.6410604034868939
15 13 15 t -15.1183272098436894      15 7 15 t 208.7485670110528013       15 3 13 t -204.0903398869704376
15 11 N 2.5654368677406967           15 5 N -2.0821819626812762           15 3 15 t 44.0416410937385712
15 11 0 B 10.9080156643548590        15 5 0 B 8.1734476146931911          15 1 N -1.9909197834262070
15 11 1 B 2.2826321158868497         15 5 1 B 2.9912119954190000          15 1 0 B 0.3003969805800217
15 11 2 B 1.0000000000000000         15 5 2 B 2.2905716060067109          15 1 1 B 0.5605904739352204
15 11 11 t 85.6713546909490324       15 5 3 B 1.5796375796515800          15 1 2 B 0.5633728436900682
15 11 13 t -135.4872079513068012     15 5 4 B 1.2413937670060894          15 1 3 B 0.7584880196839008
15 11 15 t 57.6698348943322519       15 5 5 B 1.0000000000000000          15 1 4 B 0.7724165555262611
15 9 N -2.3289382288704906           15 5 5 t 141.2270363945163140        15 1 5 B 0.8757033215822437
15 9 0 B 16.1734243381129460         15 5 7 t -447.7130171437866589       15 1 6 B 0.9351540335437449
15 9 1 B 3.6229045110143351          15 5 9 t 774.3640657773466514        15 1 7 B 1.0000000000000000
15 9 2 B 1.8379633413059982          15 5 11 t -776.1426785582560265      15 1 1 t 7.0779371051306406
15 9 3 B 1.0000000000000000          15 5 13 t 421.1587686995956359       15 1 3 t 42.9146532799192366
15 9 9 t 177.8360887940209842        15 5 15 t -95.6154155746720531       15 1 5 t -127.3672409396872028
15 9 11 t -414.0005193192334007      15 3 N 2.0185786915213594            15 1 7 t 370.8894593197135028
15 9 13 t 355.1211596278537232       15 3 0 B 2.5297284215342365          15 1 9 t -678.9190208674250806
15 9 15 t -107.9611548150770303      15 3 1 B 1.5629512291416180          15 1 11 t 726.0349793589300045
15 7 N 2.1746908116202482            15 3 2 B 1.4815640961701212          15 1 13 t -414.8296734192483416
15 7 0 B 17.1651157071880331         15 3 3 B 1.2582144623108959          15 1 15 t 97.7608510645906898
15 7 1 B 3.1903666476594628          15 3 4 B 1.1572057765542527
\end{Verbatim}
\begin{Verbatim}
16 16 N 4.5736459740053368           16 8 2 B 2.9395548029184661          16 4 5 B 1.1360455667867760
16 16 0 B 1.0000000000000000         16 8 3 B 1.5572098789564078          16 4 6 B 1.0000000000000000
16 16 16 t 1.5707963267948966        16 8 4 B 1.0000000000000000          16 4 4 t 116.1869246842475075
16 14 N -3.0710827731666776          16 8 8 t 358.3271898260139546        16 4 6 t -424.2664022234708436
16 14 0 B 4.4985542057014307         16 8 10 t -1224.2237495076239176     16 4 8 t 1032.1764534217215829
16 14 1 B 1.0000000000000000         16 8 12 t 1772.0499714559656538      16 4 10 t -1576.8505840058716884
16 14 14 t 21.1989372666106233       16 8 14 t -1216.6793096082311378     16 4 12 t 1442.2052090587189982
16 14 16 t -16.4865482862259334      16 8 16 t 324.6630647750295165       16 4 14 t -719.1633129111771497
16 12 N 2.6272529010392738           16 6 N -2.1156895380560250           16 4 16 t 150.1320642241652491
16 12 0 B 12.6522573336323330        16 6 0 B 13.8237039294767001         16 2 N -1.9999999006478862
16 12 1 B 2.4223395736994447         16 6 1 B 4.1257543581921530          16 2 0 B 1.0000000000000000
16 12 2 B 1.0000000000000000         16 6 2 B 2.8230703541621630          16 2 1 B 1.0000000000000000
16 12 12 t 99.3705967266673072       16 6 3 B 1.8090386106994828          16 2 2 B 1.0000000000000000
16 12 14 t -160.6911724071645789     16 6 4 B 1.3196822386903074          16 2 3 B 1.0000000000000000
16 12 16 t 69.1745573144717548       16 6 5 B 1.0000000000000000          16 2 4 B 1.0000000000000000
16 10 N -2.3791640253366331          16 6 6 t 238.8564569063439705        16 2 5 B 1.0000000000000000
16 10 0 B 20.7655558887621124        16 6 8 t -837.8426960208749962       16 2 6 B 1.0000000000000000
16 10 1 B 4.1082465232110560         16 6 10 t 1450.5977547062245602      16 2 7 B 1.0000000000000000
16 10 2 B 1.9536737548428559         16 6 12 t -1400.7222245380054409     16 2 2 t 23.5619449019234493
16 10 3 B 1.0000000000000000         16 6 14 t 720.7520452280333773       16 2 4 t 0.0000000000000000
16 10 10 t 228.3292123974516272      16 6 16 t -154.3625766869776081      16 2 6 t 0.0000000000000000
16 10 12 t -549.4700476795704141     16 4 N 2.0410880229411543            16 2 8 t 0.0000000000000000
16 10 14 t 478.3977528819037261      16 4 0 B 5.6897610418903820          16 2 10 t 0.0000000000000000
16 10 16 t -146.2613433122206629     16 4 1 B 2.2269869341119369          16 2 12 t 0.0000000000000000
16 8 N 2.2149276544841149            16 4 2 B 2.1339766151566397          16 2 14 t 0.0000000000000000
16 8 0 B 25.3464637800169015         16 4 3 B 1.5497521625340835          16 2 16 t 0.0000000000000000
16 8 1 B 3.6974347595021664          16 4 4 B 1.3217272019502810
\end{Verbatim}
\begin{Verbatim}
17 17 N 4.7043550333257235           17 9 15 t -1805.9115458803560354     17 3 2 B 1.6121552559691749
17 17 0 B 1.0000000000000000         17 9 17 t 479.0273641961611958       17 3 3 B 1.2975348518846115
17 17 17 t 1.5707963267948966        17 7 N -2.1503104005173516           17 3 4 B 1.2442491462219164
17 15 N -3.1488103718872846          17 7 0 B 22.2056917453553953         17 3 5 B 1.1335132804220121
17 15 0 B 4.7888694378633948         17 7 1 B 5.5742712671096645          17 3 6 B 1.0621731609209809
17 15 1 B 1.0000000000000000         17 7 2 B 3.4097299654051860          17 3 7 B 1.0000000000000000
17 15 15 t 22.5670155674884858       17 7 3 B 2.0573057581376200          17 3 3 t 63.4599501576381074
17 15 17 t -17.8546265871037959      17 7 4 B 1.3973615132742749          17 3 5 t -172.6462721800493141
17 13 N 2.6880025645920852           17 7 5 B 1.0000000000000000          17 3 7 t 500.9144593742275045
17 13 0 B 14.5250715206234363        17 7 7 t 383.6868093029841306        17 3 9 t -1065.9409055926598904
17 13 1 B 2.5619593186783039         17 7 9 t -1436.8515808135399268      17 3 11 t 1512.5387606761301151
17 13 2 B 1.0000000000000000         17 7 11 t 2499.6972737766235672      17 3 13 t -1333.2169861874255585
17 13 13 t 114.0796449551422862      17 7 13 t -2359.3735134010504545     17 3 15 t 654.8550562802460922
17 13 15 t -187.9161270385062156     17 7 15 t 1169.3508504855376702      17 3 17 t -136.4021176261836069
17 13 17 t 81.6904637173384125       17 7 17 t -239.2310797558111239      17 1 N 1.9917092553400321
17 11 N -2.4289438063123746          17 5 N 2.0659326829614874            17 1 0 B 0.2756152620524348
17 11 0 B 26.1614689110649345        17 5 0 B 11.7502669760613783         17 1 1 B 0.5390264985716238
17 11 1 B 4.6209862043571837         17 5 1 B 2.8698722718175076          17 1 2 B 0.5059278191935087
17 11 2 B 2.0691852925844865         17 5 2 B 3.0442155743550820          17 1 3 B 0.7540217889688126
17 11 3 B 1.0000000000000000         17 5 3 B 1.8584432521331634          17 1 4 B 0.6910432520651795
17 11 11 t 287.6603748834177810      17 5 4 B 1.4991298453821514          17 1 5 B 0.8487716778561311
17 11 13 t -710.5499333747960728     17 5 5 B 1.2023096017095691          17 1 6 B 0.8775544449283664
17 11 15 t 626.3743838973832351      17 5 6 B 1.0000000000000000          17 1 7 B 0.9460698537633798
17 11 17 t -192.4892511184406669     17 5 5 t 239.9445906631356754        17 1 8 B 1.0000000000000000
17 9 N 2.2555368721092498            17 5 7 t -1088.0447277893797048      17 1 1 t 7.3599025010898162
17 9 0 B 36.1397690620288935         17 5 9 t 2773.5140942227753696       17 1 3 t 56.2720947754747520
17 9 1 B 4.2233101516178884          17 5 11 t -4253.4995954180588327     17 1 5 t -221.7001828928804588
17 9 2 B 3.3708383693267482          17 5 13 t 3859.8840314946887051      17 1 7 t 863.8953199379054817
17 9 3 B 1.6468260136460510          17 5 15 t -1905.4600451874390623     17 1 9 t -2186.9593881337672724
17 9 4 B 1.0000000000000000          17 5 17 t 394.0820042626115058       17 1 11 t 3430.4617108494596760
17 9 9 t 510.9139484446574903        17 3 N -2.0160651721504764           17 1 13 t -3214.7238938868461365
17 9 11 t -1804.8332311478587898     17 3 0 B 2.6933239349209086          17 1 15 t 1645.9481169519705495
17 9 13 t 2634.9406313285502087      17 3 1 B 1.6465617397268659          17 1 17 t -353.8501425468931652
\end{Verbatim}
\begin{Verbatim}
18 18 N 4.8315367958275166           18 10 16 t -2573.1552331500227063    18 4 2 B 2.3743580769488602
18 18 0 B 1.0000000000000000         18 10 18 t 678.9350911071790683      18 4 3 B 1.6948363286471542
18 18 18 t 1.5707963267948966        18 8 N -2.1855197337991643           18 4 4 B 1.5092228845196004
18 16 N -3.2247975804116264          18 8 0 B 34.2475588385592515         18 4 5 B 1.2747270837648912
18 16 0 B 5.0791601251589606         18 8 1 B 7.3743592134050420          18 4 6 B 1.1228183248124760
18 16 1 B 1.0000000000000000         18 8 2 B 4.0568820861948581          18 4 7 B 1.0000000000000000
18 16 16 t 23.9349782034084080       18 8 3 B 2.3236947223819913          18 4 4 t 147.0070327903731425
18 16 18 t -19.2225892230237182      18 8 4 B 1.4745987561177501          18 4 6 t -589.5809340309118604
18 14 N 2.7476926263806545           18 8 5 B 1.0000000000000000          18 4 8 t 1625.1722935042603754
18 14 0 B 16.5264625230186114        18 8 8 t 591.7553358783106491        18 4 10 t -3029.7161979529974655
18 14 1 B 2.7015179076558746         18 8 10 t -2321.6777793230023732     18 4 12 t 3758.1223154666113849
18 14 2 B 1.0000000000000000         18 8 12 t 4070.1366612247420288      18 4 14 t -2960.4903362927947935
18 14 14 t 129.7985331303557717      18 8 14 t -3796.3910417231832564     18 4 16 t 1333.2439096580134046
18 14 16 t -217.1617221995467186     18 8 16 t 1837.7027247009313719      18 4 18 t -260.1961382406307389
18 14 18 t 95.2171707031654300       18 8 18 t -364.2471411630545573      18 2 N 1.9999999006447436
18 12 N -2.4783771598962308          18 6 N 2.0923013947553922            18 2 0 B 1.0000000000000000
18 12 0 B 32.3738708633392437        18 6 0 B 22.6520356499542201         18 2 1 B 1.0000000000000000
18 12 1 B 5.1585495021121128         18 6 1 B 3.2815016642368292          18 2 2 B 1.0000000000000000
18 12 2 B 2.1843016891646023         18 6 2 B 4.3155956425062988          18 2 3 B 1.0000000000000000
18 12 3 B 1.0000000000000000         18 6 3 B 2.1698894974038328          18 2 4 B 1.0000000000000000
18 12 12 t 355.9693020538592878      18 6 4 B 1.6924947646917945          18 2 5 B 1.0000000000000000
18 12 14 t -897.7442633619235235     18 6 5 B 1.2678064584702590          18 2 6 B 1.0000000000000000
18 12 16 t 799.6335750312537358      18 6 6 B 1.0000000000000000          18 2 7 B 1.0000000000000000
18 12 18 t -246.8630394356252239     18 6 6 t 462.5625471138767887        18 2 8 B 1.0000000000000000
18 10 N 2.2962757301060572           18 6 8 t -2373.3187633600057159      18 2 2 t 26.7035375555132425
18 10 0 B 50.0420726075765391        18 6 10 t 6250.0453578121631674      18 2 4 t 0.0000000000000000
18 10 1 B 4.7661539722646569         18 6 12 t -9632.0465823872124024     18 2 6 t 0.0000000000000000
18 10 2 B 3.8319466896276195         18 6 14 t 8709.0371162825689597      18 2 8 t 0.0000000000000000
18 10 3 B 1.7360757931538117         18 6 16 t -4275.5630557473350349     18 2 10 t 0.0000000000000000
18 10 4 B 1.0000000000000000         18 6 18 t 879.7037325342778936       18 2 12 t 0.0000000000000000
18 10 10 t 707.4531345346626755      18 4 N -2.0364217993890817           18 2 14 t 0.0000000000000000
18 10 12 t -2560.2928806460504714    18 4 0 B 6.2391722500960611          18 2 16 t 0.0000000000000000
18 10 14 t 3761.1970550953855032     18 4 1 B 2.6645162499529649          18 2 18 t 0.0000000000000000
\end{Verbatim}
\begin{Verbatim}
19 19 N 4.9554617130434807           19 9 4 B 1.5515177974697360          19 3 N 2.0140344461010965
19 19 0 B 1.0000000000000000         19 9 5 B 1.0000000000000000          19 3 0 B 2.8870851432304026
19 19 19 t 1.5707963267948966        19 9 9 t 882.7155869998515411        19 3 1 B 1.7327274960463565
19 17 N -3.2991432317697549          19 9 11 t -3587.6624868369859277     19 3 2 B 1.6812341170956166
19 17 0 B 5.3694306559472783         19 9 13 t 6347.9756961414618347      19 3 3 B 1.4061631611334495
19 17 1 B 1.0000000000000000         19 9 15 t -5894.5710826738308006     19 3 4 B 1.2963670269278264
19 17 17 t 25.3028458540256912       19 9 17 t 2816.3286044191924321      19 3 5 B 1.1975673618221835
19 17 19 t -20.5904568736410013      19 9 19 t -547.5075584549452169      19 3 6 B 1.1163277653165834
19 15 N 2.8063508102850436           19 7 N 2.1197410362558726            19 3 7 B 1.0537992739280148
19 15 0 B 18.6554395964246721        19 7 0 B 41.1777738936734956         19 3 8 B 1.0000000000000000
19 15 1 B 2.8410464101819032         19 7 1 B 3.1688727220988675          19 3 3 t 77.0953865482173856
19 15 2 B 1.0000000000000000         19 7 2 B 6.0585298875672012          19 3 5 t -246.6034622725846592
19 15 15 t 146.5194799640397177      19 7 3 B 2.4738534928836400          19 3 7 t 824.6105673512483028
19 15 17 t -248.4119072744038279     19 7 4 B 1.9031073996044318          19 3 9 t -1983.5586192630334240
19 15 19 t 109.7464089443385933      19 7 5 B 1.3327278100773829          19 3 11 t 3206.3096984624267550
19 13 N -2.5268604461591267          19 7 6 B 1.0000000000000001          19 3 13 t -3377.2444541902161232
19 13 0 B 39.6243788663919271        19 7 7 t 840.8646477110504918        19 3 15 t 2214.9819498829082082
19 13 1 B 5.7289705624915667         19 7 9 t -4656.9309029739341360      19 3 17 t -818.7571474163695992
19 13 2 B 2.2998105347513772         19 7 11 t 12527.4445153716065566     19 3 19 t 129.8696184529163965
19 13 3 B 1.0000000000000000         19 7 13 t -19347.4227912590045751    19 1 N -1.9932153327353564
19 13 13 t 435.6928014240043810      19 7 15 t 17416.8724476159500268     19 1 0 B 0.2500157227307733
19 13 15 t -1118.0984400385783844    19 7 17 t -8498.4993263403619658     19 1 1 B 0.5407936670952316
19 13 17 t 1004.9816885516887013     19 7 19 t 1738.0917621230272588      19 1 2 B 0.4298580547752993
19 13 19 t -311.5804756495504216     19 5 N -2.0595601051344954           19 1 3 B 0.7771209139552718
19 11 N 2.3369787616933805           19 5 0 B 12.9217242351743188         19 1 4 B 0.6047459579137593
19 11 0 B 67.5893740400360286        19 5 1 B 4.3258584076192378          19 1 5 B 0.8445103405404990
19 11 1 B 5.3243452138082515         19 5 2 B 3.2204010271515551          19 1 6 B 0.8107762405825417
19 11 2 B 4.3226135409397603         19 5 3 B 2.2524600413152776          19 1 7 B 0.9043309034267873
19 11 3 B 1.8250555024343184         19 5 4 B 1.7860868567988921          19 1 8 B 0.9481811439960454
19 11 4 B 1.0000000000000000         19 5 5 B 1.4256807590239963          19 1 9 B 1.0000000000000000
19 11 11 t 955.5222642520944193      19 5 6 B 1.1821689404862805          19 1 1 t 7.4617517992229321
19 11 13 t -3521.0044284486100400    19 5 7 B 1.0000000000000000          19 1 3 t 78.1047505034331614
19 11 15 t 5196.5367553354146831     19 5 5 t 304.4609544670262234        19 1 5 t -431.6109626162709538
19 11 17 t -3548.9448250742338173    19 5 7 t -1417.7472190922560317      19 1 7 t 2155.7916617608702697
19 11 19 t 932.0274008764888243      19 5 9 t 3706.2604135866805950       19 1 9 t -6910.8275564899017612
19 9 N -2.2209969577161435           19 5 11 t -6063.6953093415627954     19 1 11 t 14091.9029412855347747
19 9 0 B 51.0867450964692212         19 5 13 t 6363.9168808467276867      19 1 13 t -18236.8302356848543080
19 9 1 B 9.5598928109806252          19 5 15 t -4194.2266384523506016     19 1 15 t 14510.3832822368307653
19 9 2 B 4.7716481861512670          19 5 17 t 1590.4178331652397128      19 1 17 t -6474.5295352441212662
19 9 3 B 2.6074612141962221          19 5 19 t -265.8249702775813401      19 1 19 t 1239.9990326583594218
\end{Verbatim}
\begin{Verbatim}
20 20 N 5.0763674162946095           20 10 4 B 1.6282116478783229         20 4 N 2.0312252252126151
20 20 0 B 1.0000000000000000         20 10 5 B 1.0000000000000000         20 4 0 B 7.3579885093574164
20 20 20 t 1.5707963267948966        20 10 10 t 1280.3429046044858101     20 4 1 B 2.6908666893965620
20 18 N -3.3719393271319665          20 10 12 t -5351.0680117255475233    20 4 2 B 2.7404820309713013
20 18 0 B 5.6596844349598239         20 10 14 t 9561.8131750057189344     20 4 3 B 1.9297990615138486
20 18 1 B 1.0000000000000000         20 10 16 t -8879.9976498699416604    20 4 4 B 1.6340500570805511
20 18 18 t 26.6706345637594242       20 10 18 t 4217.7806867677439980     20 4 5 B 1.4201230975512907
20 18 20 t -21.9582455833747343      20 10 20 t -811.5923451877156959     20 4 6 B 1.2368386606193576
20 16 N 2.8639778109447716           20 8 N 2.1479567583995097            20 4 7 B 1.1071621717570473
20 16 0 B 20.9152873025653696        20 8 0 B 71.1404645469500162         20 4 8 B 1.0000000000000000
20 16 1 B 2.9804618321822281         20 8 1 B 2.1561360490877633          20 4 4 t 196.4843224926606691
20 16 2 B 1.0000000000000000         20 8 2 B 8.3876104207186000          20 4 6 t -997.0293023638399100
20 16 16 t 164.2682823436481206      20 8 3 B 2.7636213762468689          20 4 8 t 3526.6999021114184412
20 16 18 t -281.7195797058539260     20 8 4 B 2.1312300347020578          20 4 10 t -8339.8882666480311324
20 16 20 t 125.3052789961802885      20 8 5 B 1.3972249772447091          20 4 12 t 12995.5100044957628910
20 14 N -2.5749353453241394          20 8 6 B 1.0000000000000000          20 4 14 t -13100.6024843265992431
20 14 0 B 47.8138015893461020        20 8 8 t 1452.7133451588115117       20 4 16 t 8187.9689741312664230
20 14 1 B 6.3240219163850612         20 8 10 t -8452.1057252466537440     20 4 18 t -2873.0768830200843299
20 14 2 B 2.4149565879691125         20 8 12 t 23038.9978385441140418     20 4 20 t 430.6372706829594336
20 14 3 B 1.0000000000000000         20 8 14 t -35560.5185652064349126    20 2 N -1.9999999006412079
20 14 14 t 525.7402073465139288      20 8 16 t 31830.7327599016509718     20 2 0 B 1.0000000000000000
20 14 16 t -1368.6118637061521686    20 8 18 t -15420.4616658470120994    20 2 1 B 1.0000000000000000
20 14 18 t 1239.6646090655339399     20 8 20 t 3131.0623649438578866      20 2 2 B 1.0000000000000000
20 14 20 t -385.7973784183314237     20 6 N -2.0846001565800639           20 2 3 B 1.0000000000000000
20 12 N 2.3770446643969038           20 6 0 B 24.4939312811618195         20 2 4 B 1.0000000000000000
20 12 0 B 89.9248696271832783        20 6 1 B 7.1105657139718798          20 2 5 B 1.0000000000000000
20 12 1 B 5.8301988069098267         20 6 2 B 4.0128635606680140          20 2 6 B 1.0000000000000000
20 12 2 B 4.8475170778200022         20 6 3 B 3.0525709862291979          20 2 7 B 1.0000000000000000
20 12 3 B 1.9134202169672362         20 6 4 B 2.0605205299564764          20 2 8 B 1.0000000000000000
20 12 4 B 1.0000000000000000         20 6 5 B 1.5891193618719644          20 2 9 B 1.0000000000000000
20 12 12 t 1271.2828940810051238     20 6 6 B 1.2403708793813411          20 2 2 t 29.8451302091030358
20 12 14 t -4755.4416009904144911    20 6 7 B 1.0000000000000000          20 2 4 t 0.0000000000000000
20 12 16 t 7049.8083875606330179     20 6 6 t 577.1246592782340339        20 2 6 t 0.0000000000000000
20 12 18 t -4810.2221843306628251    20 6 8 t -2867.1013119288584085      20 2 8 t 0.0000000000000000
20 12 20 t 1258.7096706205932442     20 6 10 t 7068.5582991535536793      20 2 10 t 0.0000000000000000
20 10 N -2.2565354917881731          20 6 12 t -10018.2770630370964921    20 2 12 t 0.0000000000000000
20 10 0 B 74.0992371347050610        20 6 14 t 8229.4339163159072562      20 2 14 t 0.0000000000000000
20 10 1 B 12.1611335065566218        20 6 16 t -3591.0229283835716819     20 2 16 t 0.0000000000000000
20 10 2 B 5.5615692138034645         20 6 18 t 591.3749730510858504       20 2 18 t 0.0000000000000000
20 10 3 B 2.9079798702195933         20 6 20 t 33.4714004526692119        20 2 20 t 0.0000000000000000
\end{Verbatim}
\begin{Verbatim}
21 21 N 5.1944640246384188           21 9 N 2.1767364830900448            21 5 17 t 30678.6796211724084849
21 21 0 B 1.0000000000000000         21 9 0 B 117.5828226938516972        21 5 19 t -11397.8933038543720282
21 21 21 t 1.5707963267948966        21 9 1 B -0.2173839640798862         21 5 21 t 1836.2884368396630813
21 19 N -3.4432711347181069          21 9 2 B 11.4211336476051441         21 3 N -2.0125684392941655
21 19 0 B 5.9499241429846107         21 9 3 B 3.0348469555180011          21 3 0 B 3.0392947560476160
21 19 1 B 1.0000000000000000         21 9 4 B 2.3766409958773836          21 3 1 B 1.8025805214506147
21 19 19 t 28.0383569655254993       21 9 5 B 1.4614144695776219          21 3 2 B 1.8141584769952755
21 19 21 t -23.3259679851408095      21 9 6 B 1.0000000000000000          21 3 3 B 1.4115146247840909
21 17 N 2.9206343954976344           21 9 9 t 2401.0826577618121398       21 3 4 B 1.4003773130598356
21 17 0 B 23.3023565455075667        21 9 11 t -14433.1302892887751557    21 3 5 B 1.2416875761267185
21 17 1 B 3.1198638711025562         21 9 13 t 39647.7651624075868008     21 3 6 B 1.1777179354727603
21 17 2 B 1.0000000000000000         21 9 15 t -61041.9580349496271329    21 3 7 B 1.1047041285080724
21 17 17 t 183.0162803367415100      21 9 17 t 54282.3238514186300901     21 3 8 B 1.0489047028245642
21 17 19 t -317.0258535852029997     21 9 19 t -26091.6231582482026781    21 3 9 B 1.0000000000000000
21 17 21 t 141.8635548824359728      21 9 21 t 5255.9601631469095918      21 3 3 t 90.7081477380851463
21 15 N -2.6224382096571024          21 7 N -2.1108973753148865           21 3 5 t -332.1890762669883359
21 15 0 B 57.0330561498071093        21 7 0 B 43.3098159854721585         21 3 7 t 1341.1959463362352569
21 15 1 B 6.9511403469125969         21 7 1 B 11.7408870163192491         21 3 9 t -4167.9075263034353442
21 15 2 B 2.5300482563144435         21 7 2 B 4.5542918154707938          21 3 11 t 9281.7917902305145995
21 15 3 B 1.0000000000000000         21 7 3 B 4.1853330385215114          21 3 13 t -14338.8462037243627672
21 15 15 t 627.1112057420286748      21 7 4 B 2.3197143406573213          21 3 15 t 14889.3913690825377316
21 15 17 t -1652.0382770227977007    21 7 5 B 1.7667980059396740          21 3 17 t -9876.7945607285286259
21 15 19 t 1506.2009374797931560     21 7 6 B 1.2976026173901475          21 3 19 t 3773.7959688090845459
21 15 21 t -470.2782919114598537     21 7 7 B 1.0000000000000000          21 3 21 t -631.3007249640391702
21 13 N 2.4178214212084064           21 7 7 t 1020.4634979621384335       21 1 N 1.9937657954504884
21 13 0 B 116.0000455462034282       21 7 9 t -5206.7775548874111853      21 1 0 B 0.2254326379594920
21 13 1 B 6.4769147846058111         21 7 11 t 12064.3993014110389838     21 1 1 B 0.5447529107807667
21 13 2 B 5.3920911472965315         21 7 13 t -14485.0450866743426438    21 1 2 B 0.3461845753702508
21 13 3 B 2.0024393338529108         21 7 15 t 7628.5103535572913039      21 1 3 B 0.8436019369945777
21 13 4 B 1.0000000000000000         21 7 17 t 926.8894540144079766       21 1 4 B 0.4554123811693986
21 13 13 t 1639.9120090681534704     21 7 19 t -2809.6869958402479690     21 1 5 B 0.9047041641432682
21 13 15 t -6193.3871339782086665    21 7 21 t 884.8089753590485496       21 1 6 B 0.7148960816774482
21 13 17 t 9198.0627037932054412     21 5 N 2.0501204297564524            21 1 7 B 0.8927150756374579
21 13 19 t -6262.3767653155509760    21 5 0 B 17.3705295776085946         21 1 8 B 0.8961740320332298
21 13 21 t 1631.9263533735548005     21 5 1 B 3.4096958316608328          21 1 9 B 0.9557635425418766
21 11 N -2.2920144019731056          21 5 2 B 4.6422039246361369          21 1 10 B 1.0000000000000000
21 11 0 B 104.8966171171646979       21 5 3 B 2.4186838504293042          21 1 1 t 7.4362839525755294
21 11 1 B 15.2094382333456622        21 5 4 B 2.0920309445948533          21 1 3 t 105.3332934399485422
21 11 2 B 6.4335886519623850         21 5 5 B 1.6479739909417526          21 1 5 t -768.7551597069500094
21 11 3 B 3.2248110540500107         21 5 6 B 1.3664291200792726          21 1 7 t 4805.0082357569079016
21 11 4 B 1.7047413112687385         21 5 7 B 1.1597847898422729          21 1 9 t -19364.7914362667575615
21 11 5 B 1.0000000000000000         21 5 8 B 1.0000000000000000          21 1 11 t 50708.4189157123040925
21 11 11 t 1812.4834294693828387     21 5 5 t 463.8545889348246875        21 1 13 t -87314.0281364368200706
21 11 13 t -7748.4160133214866630    21 5 7 t -2982.4291859295414455      21 1 15 t 98136.1314582401562442
21 11 15 t 13980.4740750794212030    21 5 9 t 11360.0472829041162358      21 1 17 t -69312.9199906117407152
21 11 17 t -13018.5530905891450990   21 5 11 t -27888.2325313027807835    21 1 19 t 27928.1084091622145907
21 11 19 t 6174.2099383222637950     21 5 13 t 44866.1775777882984842     21 1 21 t -4896.9551503791457150
21 11 21 t -1182.9195793656922120    21 5 15 t -46909.7889489971034740
\end{Verbatim}
\begin{Verbatim}
22 22 N 5.3099383998567259           22 10 N 2.2059098952360651           22 6 18 t 90622.2938729998348969
22 22 0 B 1.0000000000000000         22 10 0 B 187.0253584343142978       22 6 20 t -34594.0769018235990415
22 22 22 t 1.5707963267948966        22 10 1 B -4.5086360950784801        22 6 22 t 5743.0390458716492482
22 20 N -3.5132175061024411          22 10 2 B 15.2837130062012351        22 4 N -2.0289772625800599
22 20 0 B 6.2401519189383327         22 10 3 B 3.2842327100106702         22 4 0 B 7.7004403281345169
22 20 1 B 1.0000000000000000         22 10 4 B 2.6390023250820694         22 4 1 B 3.3534260837572330
22 20 20 t 29.4060231387313756       22 10 5 B 1.5253820323431407         22 4 2 B 2.7849793834844465
22 20 22 t -24.6936341583466857      22 10 6 B 1.0000000000000000         22 4 3 B 2.1790687628271992
22 18 N 2.9763372661652544           22 10 10 t 3819.1236985995578729     22 4 4 B 1.7921207857350834
22 18 0 B 25.8188688106929741        22 10 12 t -23467.1498149236719385   22 4 5 B 1.5598714124203373
22 18 1 B 3.2592152795617161         22 10 14 t 64730.3756443610075748    22 4 6 B 1.3618967149233884
22 18 2 B 1.0000000000000000         22 10 16 t -99291.1746241201218646   22 4 7 B 1.2153164679631206
22 18 18 t 202.7809214491792235      22 10 18 t 87684.2530323728821561    22 4 8 B 1.0970045955458117
22 18 20 t -354.3662090046649889     22 10 20 t -41808.5848563260604739   22 4 9 B 1.0000000000000000
22 18 22 t 159.4392691894602486      22 10 22 t 8353.5772722847403301     22 4 4 t 229.8206442606027635
22 16 N -2.6691043243520512          22 8 N -2.1380298704137577           22 4 6 t -1167.6348552983911421
22 16 0 B 67.4930223186078415        22 8 0 B 72.4441550048697167         22 4 8 t 4059.7862528768185031
22 16 1 B 7.5959273906417015         22 8 1 B 19.1498238021675630         22 4 10 t -9566.7562989353965967
22 16 2 B 2.6451790531269457         22 8 2 B 4.6260416512381420          22 4 12 t 15314.4230563618011239
22 16 3 B 1.0000000000000000         22 8 3 B 5.7349875646354725          22 4 14 t -16520.3757123593265853
22 16 16 t 742.1245407964862189      22 8 4 B 2.5561913581394977          22 4 16 t 11676.7599148698604505
22 16 18 t -1975.8088706692233584    22 8 5 B 1.9591946243834276          22 4 18 t -5080.4240891420650256
22 16 20 t 1812.4999072966868491     22 8 6 B 1.3540974837810120          22 4 20 t 1189.3452737656462190
22 16 22 t -567.8200031363854333     22 8 7 B 1.0000000000000000          22 4 22 t -105.0990561904466744
22 14 N 2.4578032457578116           22 8 8 t 1706.9251886911421558       22 2 N 1.9999999006372789
22 14 0 B 148.1811923327921448       22 8 10 t -8790.0266676804947471     22 2 0 B 1.0000000000000000
22 14 1 B 7.0654313922330367         22 8 12 t 19183.7003520789808107     22 2 1 B 1.0000000000000000
22 14 2 B 5.9707320888796346         22 8 14 t -19081.0222139869581127    22 2 2 B 1.0000000000000000
22 14 3 B 2.0909223890697953         22 8 16 t 2653.2464824518938505      22 2 3 B 1.0000000000000000
22 14 4 B 1.0000000000000000         22 8 18 t 11663.7725200989004375     22 2 4 B 1.0000000000000000
22 14 14 t 2094.8622535479419932     22 8 20 t -9866.6775146228109393     22 2 5 B 1.0000000000000000
22 14 16 t -7979.9082817786917376    22 8 22 t 2553.6437978712699937      22 2 6 B 1.0000000000000000
22 14 18 t 11877.0067418568049660    22 6 N 2.0702342247584646            22 2 7 B 1.0000000000000000
22 14 20 t -8075.4987770681971048    22 6 0 B 38.3367306899896766         22 2 8 B 1.0000000000000000
22 14 22 t 2097.6752303832959529     22 6 1 B 2.8630598325912801          22 2 9 B 1.0000000000000000
22 12 N -2.3268892698014896          22 6 2 B 8.1139057134838064          22 2 10 B 1.0000000000000000
22 12 0 B 146.3369691384533503       22 6 3 B 2.6615801127908681          22 2 2 t 32.9867228626928290
22 12 1 B 18.7165309711903230        22 6 4 B 2.7509467628750949          22 2 4 t 0.0000000000000000
22 12 2 B 7.3972485230686150         22 6 5 B 1.8675350722616150          22 2 6 t 0.0000000000000000
22 12 3 B 3.5562839945205454         22 6 6 B 1.5072170782531688          22 2 8 t 0.0000000000000000
22 12 4 B 1.7811822083450392         22 6 7 B 1.2116140804019752          22 2 10 t 0.0000000000000000
22 12 5 B 1.0000000000000000         22 6 8 B 1.0000000000000000          22 2 12 t 0.0000000000000000
22 12 12 t 2528.5213095667873698     22 6 6 t 1023.7263277357364314       22 2 14 t 0.0000000000000000
22 12 14 t -11025.6143523400595386   22 6 8 t -7578.1800157796297998      22 2 16 t 0.0000000000000000
22 12 16 t 20095.3971026191314563    22 6 10 t 30449.6825402302207686     22 2 18 t 0.0000000000000000
22 12 18 t -18803.2365275653517911   22 6 12 t -76904.8674112709243707    22 2 20 t 0.0000000000000000
22 12 20 t 8934.0157046617550576     22 6 14 t 126496.7572773824906041    22 2 22 t 0.0000000000000000
22 12 22 t -1711.8044773475186913    22 6 16 t -135231.6711977902654936
\end{Verbatim}
\normalsize

\section{Bernstein Polynomials} \label{app.Bernst}
The Bernstein Polynomials are defined over the unit interval $0\le x \le 1$
as \cite{FaroukiCAGD29,RababahCMAM3}:
\begin{eqnarray}
B_{i,j}(x) & \equiv & \binom{j}{i} x^i (1-x)^{j-i}
=\binom{j}{i} \sum_{k=0}^{j-i}\binom{j-i}{k}(-)^k x^{i+k}
=\sum_{k=0}^{j-i}\binom{j}{k}\binom{j-k}{i}(-)^k x^{i+k} \nonumber \\
& = &\binom{j}{i} \sum_{l=i}^j\binom{j-i}{l-i}(-)^{l-i} x^l
=\sum_{l=i}^j \binom{j}{l} \binom{l}{i} (-)^{l-i} x^l
.
\end{eqnarray}
Starting from the second term of the right hand side
\begin{equation}
\frac{1}{\binom{j}{i}}B_{i,j}(x)
=\sum_{k=0}^{j-i}\binom{j-i}{k}(-)^k x^{i+k};
\end{equation}
and setting $j-i-k\equiv l$ and then $j-i\equiv n$ yields
\begin{equation}
\frac{1}{\binom{j}{n}}B_{j-n,j}(x)
=\sum_{k=0}^n\binom{n}{k}(-)^{n+k} x^{j-k}
.
\end{equation}
The well-known inversion of the binomial series expands the monomials
in finite series of Bernstein Polynomials \cite[Table 2.1]{Riordan}\cite{RiordanAMM71,GarloffSAMV}:
\begin{equation}
x^{j-n}
=
\sum_{k=0}^n \binom{n}{k}
\frac{1}{\binom{j}{k}}B_{j-k,j}(x)
.
\end{equation}
\begin{equation}
x^i
=
\sum_{k=0}^{j-i} \binom{j-i}{k}
\frac{1}{\binom{j}{k}}B_{j-k,j}(x)
=
\frac{1}{\binom{j}{i}}
\sum_{k=0}^{j-i} \binom{j-k}{i} B_{j-k,j}(x)
=\frac{1}{\binom{j}{i}} \sum_{s=i}^{j} \binom{s}{i} B_{s,j}(x)
.
\end{equation}

\bibliographystyle{apsrmp4-1}
\bibliography{all}

\end{document}